\documentclass{LMCS}

\def\doi{8 (1:10) 2012}
\lmcsheading%
{\doi}
{1--28}
{}
{}
{Aug.~24, 2010}
{Feb.~27, 2012}
{}

\usepackage[latin1]{inputenc}

\usepackage{amssymb}
\usepackage{amsthm}

\usepackage{logique_lmcs}
\usepackage{enumerate,hyperref}

\newcommand{\NN}{\mathbb{N}}

\newcommand{\nn}{\mbox{\sffamily n}}

\newcommand{\inde}{^{\mbox{\footnotesize ent}}}
\newcommand{\ppt}{{\scriptstyle<}}

\newcommand{\et}{{\scriptstyle\land}}
\newcommand{\ou}{{\scriptstyle\lor}}
\newcommand{\non}{{\scriptstyle\neg}}

\newcommand{\ZFe}{ZF$_\varepsilon$}
\newcommand{\gd}{\gimel\mathbf{2}}
\newcommand{\gn}{\gimel\mathbb{N}}
\newcommand{\mb}[1]{\mathbf{#1}}

\theoremstyle{plain}
\newtheorem{theorem}{Theorem}[section]
\newtheorem{lemma}[theorem]{Lemma}
\newtheorem{corollary}[theorem]{Corollary}
\newtheorem{proposition}[theorem]{Proposition}

\newcommand{\mapto}{\,\hookrightarrow}

\begin{document}

\title{Realizability algebras II~:\\
new models of \ ZF + DC}

\author[J.-L.~Krivine]{Jean-Louis Krivine}
\address{Université Paris VII, C.N.R.S.}
\email{krivine@pps.univ-paris-diderot.fr}


\begin{abstract}\noindent
Using the proof-program (Curry-Howard) correspondence, we give a new method to obtain models of ZF and relative consistency results in set theory. We show the relative consistency of ZF + DC + there exists a sequence of
subsets of R the cardinals of which are strictly decreasing + other similar properties of R.
These results seem not to have been previously obtained by forcing. 
\end{abstract}

\maketitle\noindent
\section*{Introduction}
The technology of \emph{classical realizability} was developed in \cite{krivine3,krivine6} in order to
extend the proof-program correspondence (also known as \emph{Curry-Howard correspondence}) from pure intuitionistic
logic to the whole of mathematical proofs, with excluded middle, axioms of~ZF\/, dependent choice, existence of a
well ordering on~${\mathcal P}(\NN)$,~\ldots\\
We show here that this technology is also a new method in order to build models of ZF and
to obtain relative consistency results.

\smallskip\noindent
The main tools are~:

\smallskip\noindent
$\bullet$~~The notion of \emph{realizability algebra} \cite{krivine6}, which comes from combinatory
logic~\cite{curry} and plays a role similar to a set of forcing conditions. The extension from intuitionistic
to classical logic was made possible by Griffin's discovery \cite{griffin} of the relation between the law of
Peirce and the instruction {\tt call-with-current-continuation} of the programming language SCHEME.\\
In this paper, we only use the simplest case of realizability algebra, which I call \emph{standard
realizability algebra}~; somewhat like the \emph{binary tree} in the case of forcing.

\smallskip\noindent
$\bullet$~~The theory \ZFe\ \cite{krivine1} which is a conservative extension of ZF\/, with a notion of
\emph{strong membership}, denoted as \ $\varepsilon$.

\smallskip\noindent
The theory \ZFe\ is essentially ZF without the extensionality  axiom. We note an analogy with the
Fraenkel-Mostowski models with ``urelements''~: we obtain a non well orderable set, which is a
Boolean algebra denoted $\gd$, all elements of which (except $1$) are empty\/. But we also notice
two important differences~:\\
$\bullet$~~The final model of ZF + $\neg$ AC is obtained directly, without taking a suitable submodel.\\
$\bullet$~~There exists an injection from the ``pathological set'' $\gd$ into $\mathbb{R}$, and
therefore \emph{$\mathbb{R}$ is also not well orderable}.

\smallskip\noindent
We show the consistency, relatively to the consistency of ZF\/, of the theory ZF + DC (dependent choice)
with the following properties~:

\smallskip
there exists a sequence $({\mathcal X}_n)_{n\in\NN}$ of infinite subsets of $\mathbb{R}$, the ``cardinals'' of which are strictly increasing (this means that there is an injection but no surjection from ${\mathcal X}_n$ to ${\mathcal X}_{n+1}$), and such that \
${\mathcal X}_m\fois{\mathcal X}_n$ is equipotent with~${\mathcal X}_{mn}$ for $m,n\ge2$~;

there exists a sequence of infinite subsets of $\mathbb{R}$, the ``cardinals'' of which are strictly decreasing.

\smallskip\noindent
More detailed properties of $\mathbb{R}$ in this model are given in theorems~\ref{X_m_fois_X_n} and~\ref{subsets_P(N)}.

\smallskip\noindent
As far as I know, these consistency results are new, and it seems they cannot be obtained by forcing.
But, in any case, the fact that the simplest non trivial realizability model (which I call the
\emph{model of threads}) has a real line with such unusual properties, is of interest in itself.
Another aspect of these results, which is interesting from the point of view of computer science,
is the following~: in~\cite{krivine6}, we introduce {\em read} and {\em write} instructions in a global memory, in order to realize a weak form of the axiom of choice (well ordering of~$\mathbb{R}$).
Therefore, what we show here, is that these instructions are \emph{indispensable}~: without them,
we can build a realizability model in which $\mathbb{R}$ is not well ordered.

\section{Standard realizability algebras}\noindent
The structure of \emph{realizability algebra}, and the particular case of \emph{standard realizability algebra} are defined in \cite{krivine6}. They are variants of the usual notion of
\emph{combinatory algebra}.
Here, we only need the \emph{standard} realizability algebras, the definition of which we recall below~:

\smallskip\noindent 
We have a countable set $\Pi_0$ which is the set of \emph{stack constants}.\\
We define recursively two sets~: $\Lbd$ (the set of \emph{terms}) and $\Pi$ (the set of \emph{stacks}).
Terms and stacks are finite sequences of elements of the set~:\\
\centerline{$\Pi_0\cup\{B,C,E,I,K,W,\Ccc,\vsig,\kk,(,),[,],\ps\}$}

\noindent
which are obtained by the following rules~:

\smallskip\noindent
$\bullet$~~$B,C,E,I,K,W,\Ccc,\vsig$ are terms (\emph{elementary combinators})~;\\
$\bullet$~~each element of $\Pi_0$ is a stack (\emph{empty stacks})~;\\
$\bullet$~~if $\xi,\eta$ are terms, then $(\xi)\eta$ is a term (this operation is called
\emph{application})~;\\
$\bullet$~~if $\xi$ is a term and $\pi$ a stack, then $\xi\ps\pi$ is a stack (this operation is called
\emph{push})~;\\
$\bullet$~~if $\pi$ is a stack, then $\kk[\pi]$ is a term.\\
A term of the form $\kk[\pi]$ is called a \emph{continuation}. From now on, it will be denoted as
$\kk_\pi$.

\smallskip\noindent
A term which does not contain any continuation (i.e. in which the symbol $\kk$ does not appear)
is called \emph{proof-like}.

\smallskip\noindent
Every stack has the form \ $\pi=\xi_1\ps\ldots\ps\xi_n\ps\pi_0$,
where $\xi_1,\ldots,\xi_n\in\Lbd$ and $\pi_0\in\Pi_0$, i.e. $\pi_0$ is a stack constant.

\smallskip\noindent
If $\xi\in\Lbd$ and $\pi\in\Pi$, the ordered pair $(\xi,\pi)$ is called a \emph{process} and
denoted as $\xi\star\pi$~;\\
$\xi$ and $\pi$ are called respectively the \emph{head} and the \emph{stack} of the process $\xi\star\pi$.\\
The set of processes $\Lbd\fois\Pi$ will also be written $\Lbd\star\Pi$.

\smallskip\noindent
{\bfseries Notation.}\\ For sake of brevity, the term $(\ldots(((\xi)\eta_1)\eta_2)\ldots)\eta_n$ will be also denoted as $(\xi)\eta_1\eta_2\ldots\eta_n$ or $\xi\eta_1\eta_2\ldots\eta_n$, if the meaning is clear.
For example~: \ $\xi\eta\zeta=(\xi)\eta\zeta=(\xi\eta)\zeta=((\xi)\eta)\zeta$.

\smallskip\noindent
We now choose a recursive bijection from $\Lbd$ onto $\NN$, which is written \ $\xi\longmapsto\nn_\xi$.\\
We put \ $\sig=(BW)(B)B$ (the characteristic property of $\sig$ is given below).\\
For each $n\in\NN$, we define \ $\ul{n}\in\Lbd$ \ recursively\/, by putting~: \ $\ul{0}=KI$~; \
$\ul{n+1}=(\sig)\ul{n}$~;\\
$\ul{n}$ is the \emph{$n$-th integer} and $\sig$ is the \emph{successor} in combinatory logic.

\smallskip\noindent
We define a preorder relation \ $\succ$ \ on $\Lbd\star\Pi$. It is the least reflexive and transitive
relation such that, for all $\xi,\eta,\zeta\in\Lbd$ and $\pi,\varpi\in\Pi$, we have~:

\smallskip\noindent
$(\xi)\eta\star\pi\succ\xi\star\eta\ps\pi$.\\
$I\star\xi\ps\pi\succ\xi\star\pi$.\\
$K\star\xi\ps\eta\ps\pi\succ\xi\star\pi$.\\
$E\star\xi\ps\eta\ps\pi\succ(\xi)\eta\star\pi$.\\
$W\star\xi\ps\eta\ps\pi\succ\xi\star\eta\ps\eta\ps\pi$.\\
$C\star\xi\ps\eta\ps\zeta\ps\pi\succ\xi\star\zeta\ps\eta\ps\pi$.\\
$B\star\xi\ps\eta\ps\zeta\ps\pi\succ(\xi)(\eta)\zeta\star\pi$.\\
$\Ccc\star\xi\ps\pi\succ\xi\star\kk_\pi\ps\pi$.\\
$\kk_\pi\star\xi\ps\varpi\succ\xi\star\pi$.\\
$\vsig\star\xi\ps\eta\ps\pi\succ\xi\star\ul{\nn}\,_\eta\ps\pi$.

\smallskip\noindent
For instance, with the definition of \ $\ul{0}$ and $\sig$ given above, we have~:\\ $\ul{0}\star\xi\ps\eta\ps\pi\succ\eta\star\pi$~; \
$\sigma\star\xi\ps\eta\ps\zeta\ps\pi\succ(\xi\eta)(\eta)\zeta\star\pi$.

\smallskip\noindent
Finally, we have a subset $\bbot$ of $\Lbd\star\Pi$ which is a final segment for this
preorder, which means that~: \ \ $\xi\star\pi\in\bbot$, \ $\xi'\star\pi'\succ\xi\star\pi$ \
$\Fl$ \ $\xi'\star\pi'\in\bbot$.\\
In other words, we ask that $\bbot$ has the following properties~:

\smallskip\noindent
$(\xi)\eta\star\pi\notin\bbot\Fl\xi\star\eta\ps\pi\notin\bbot$.\\
$I\star\xi\ps\pi\notin\bbot\Fl\xi\star\pi\notin\bbot$.\\
$K\star\xi\ps\eta\ps\pi\notin\bbot\Fl\xi\star\pi\notin\bbot$.\\
$E\star\xi\ps\eta\ps\pi\notin\bbot\Fl(\xi)\eta\star\pi\notin\bbot$.\\
$W\star\xi\ps\eta\ps\pi\notin\bbot\Fl\xi\star\eta\ps\eta\ps\pi\notin\bbot$.\\
$C\star\xi\ps\eta\ps\zeta\ps\pi\notin\bbot\Fl\xi\star\zeta\ps\eta\ps\pi\notin\bbot$.\\
$B\star\xi\ps\eta\ps\zeta\ps\pi\notin\bbot\Fl(\xi)(\eta)\zeta\star\pi\notin\bbot$.\\
$\Ccc\star\xi\ps\pi\notin\bbot\Fl\xi\star\kk_\pi\ps\pi\notin\bbot$.\\
$\kk_\pi\star\xi\ps\varpi\notin\bbot\Fl\xi\star\pi\notin\bbot$.\\
$\vsig\star\xi\ps\eta\ps\pi\notin\bbot\Fl\xi\star\ul{\nn}\,_\eta\ps\pi\notin\bbot$.

\smallskip\noindent
{\small{\bfseries Remark.} Thus, the only arbitrary elements in a standard realizability algebra
are the set $\Pi_0$ of stack constants and the set $\bbot$ of processes.}

\subsection*{\Cc-terms and $\lbd$-terms}\ \\
We call \ \emph{\Cc-term} \ a term which is built with variables, the elementary combinators
$B$, $C$, $E$, $I$, $K$, $W$, $\Ccc$, $\vsig$ and the application (binary function). A closed \Cc-term
is exactly what we have called a proof-like term.

\smallskip\noindent
Given a \Cc-term $t$ and a variable $x$, we define inductively on $t$, a new \Cc-term denoted by $\lbd x\,t$,
which does not contain $x$. To this aim, we apply the first possible case in the following list~:

\smallskip\noindent
1.~$\lbd x\,t=(K)t$\label{def_lbd} if $t$ does not contain $x$.\\
2.~$\lbd x\,x=I$.\\
3.~$\lbd x\,tu=(C\lbd x(E)t)u$ if $u$ does not contain $x$.\\
4.~$\lbd x\,tx=(E)t$ if $t$ does not contain $x$.\\
5.~$\lbd x\,tx=(W)\lbd x(E)t$ (if $t$ contains $x$).\\
6.~$\lbd x(t)(u)v=\lbd x(B)tuv$ (if $uv$ contains $x$).

\smallskip\noindent
In \cite{krivine6}, it is shown that this definition is correct. This allows us to translate
every $\lbd$-term into a $\Cc$-term. In the following, almost every $\Cc$-term will be written
as a $\lbd$-term. The fundamental property of this translation is given by
theorem~\ref{beta_red_gauche}, which is proved in \cite{krivine6}~:

\begin{theorem}\label{beta_red_gauche}
Let $t$ be a \ \Cc-term with the only variables $x_1,\ldots,x_n$~;
let \ $\xi_1,\ldots,\xi_n\in\Lbd$ and $\pi\in\Pi$. Then
$\lbd x_1\ldots\lbd x_n\,t\star\xi_1\ps\ldots\ps\xi_n\ps\pi
\succ t[\xi_1/x_1,\ldots,\xi_n/x_n]\star\pi$.
\end{theorem}\noindent

\smallskip\noindent
{\small{\bfseries Remark.} The property we need for the term $\sig$ (the \emph{successor}) is \
$\sigma\star\xi\ps\eta\ps\zeta\ps\pi\succ(\xi\eta)(\eta)\zeta\star\pi$ (to prove theorem~\ref{rec_eps}).
Therefore, by theorem~\ref{beta_red_gauche}, we could define \ $\sig=\lbd n\lbd f\lbd x(nf)(f)x$.
The definition we chose is much simpler.}

\section{The formal system}
We write formulas and proofs in the language of first order logic. This formal
language consists of~:

\smallskip\noindent
$\bullet$~~\emph{individual variables} \ $x,y,\ldots$~;\\
$\bullet$~~\emph{function symbols} $f,g,\ldots$~; each one has an arity, which is an integer~; function symbols of arity~$0$ are called \emph{constant symbols}.\\
$\bullet$~~\emph{relation symbols}~; each one has an arity~; relation symbols of arity~$0$ are called \emph{propositional constants}. We have two particular propositional constants $\top,\bot$ and three particular binary relation symbols \ $\neps,\notin,\subseteq$.

\smallskip\noindent
The \emph{terms} are built in the usual way with individual variables and function symbols.\\
{\small{\bfseries Remark.} We use the word ``term'' with two different meanings~: here as a term in a first order
language, and previously as an element of the set $\Lbd$ of a realizability algebra. I think that, with the help
of the context, no confusion is possible.}

\smallskip\noindent
The \emph{atomic formulas} are the expressions $R(t_1,\ldots,t_n)$, where $R$ is a
$n$-ary relation symbol, and $t_1,\ldots,t_n$ are terms.

\smallskip\noindent
\emph{Formulas} are built as usual, from atomic formulas, \emph{with the only logical symbols} \ $\to,\pt$~:\\
$\bullet$~~each atomic formula is a formula~;\\
$\bullet$~~if $A,B$ are formulas, then $A\to B$ is a formula~;\\
$\bullet$~~if $A$ is a formula and $x$ an individual variable, then $\pt x\,A$ is a
formula.

\smallskip\noindent
{\bfseries Notations.}\\
The formula $A_1\to(A_2\to(\cdots(A_n\to B)\cdots))$ will be written
$A_1,A_2,\ldots,A_n\to B$.\\
The usual logical symbols are defined as follows~:\\
$\neg A\equiv A\to\bot$~; $A\lor B\equiv(A\to\bot),(B\to\bot)\to\bot$~;
$A\land B\equiv(A,B\to\bot)\to\bot$~;
$\ex x\,F\equiv\pt x\,(F\to\bot)\to\bot$.\\
More generally, we shall write \ $\ex x\{F_1,\ldots,F_k\}$ \ for \
$\pt x(F_1,\ldots,F_k\to\bot)\to\bot$.\\
We shall sometimes write \ $\vec{F}$ \ for a finite sequence of formulas \ $F_1,\ldots,F_k$~;\\
Then, we shall also write \ $\vec{F}\to G$ \ for  \ $F_1,\ldots,F_k\to G$ \ and \
$\ex x\{\vec{F}\}$ \ for \
$\pt x(\vec{F}\to\bot)\to\bot$.\\
$A\dbfl B$ \ is the pair of formulas \ $\{A\to B,B\to A\}$.

\smallskip\noindent
The rules of natural deduction are the following (the $A_i$'s are formulas, the $x_i$'s
are variables of \Cc-term, $t,u$ are \Cc-terms, written as $\lbd$-terms)~:

\smallskip\noindent
1.~$x_1:A_1,\ldots,x_n:A_n\vdash x_i:A_i$.

\noindent
2.~$x_1:A_1,\ldots,x_n:A_n\vdash t:A\to B$, \ \ $x_1:A_1,\ldots,x_n:A_n\vdash u:A$\\
\hspace*{\fill}$\Fl$ \ \ $x_1:A_1,\ldots,x_n:A_n\vdash tu:B$.

\noindent
3.~$x_1:A_1,\ldots,x_n:A_n,x:A\vdash t:B$ \ \ $\Fl$ \ \
$x_1:A_1,\ldots,x_n:A_n\vdash\lbd x\,t:A\to B$.\\
4.~$x_1:A_1,\ldots,x_n:A_n\vdash t:A$ \ \ $\Fl$ \ \
$x_1:A_1,\ldots,x_n:A_n\vdash t:\pt x\,A$ \ \
where $x$ is an individual variable which does not appear in $A_1,\ldots,A_n$.\\
5.~$x_1:A_1,\ldots,x_n:A_n\vdash t:\pt x\,A$ \ \ $\Fl$ \ \
$x_1:A_1,\ldots,x_n:A_n\vdash t:A[\tau/x]$ \ \ where $x$ is an individual variable and
$\tau$ is a term.\\
6.~$x_1:A_1,\ldots,x_n:A_n\vdash\Ccc:((A\to B)\to A)\to A$ \ (law of Peirce).\\
7.~$x_1:A_1,\ldots,x_n:A_n\vdash t:\bot$ \ \ $\Fl$ \ \
$x_1:A_1,\ldots,x_n:A_n\vdash t:A$ \ \ for every formula $A$.

\section{The theory \ZFe}
We write below a set of axioms for a theory called \ \ZFe. Then~:\\
$\bullet$ We show that \ZFe\ is a conservative extension of ZF.\\
$\bullet$ We define the \emph{realizability models} and we show that each axiom of
\ZFe\ is realized by a proof-like $\Cc$-term, in every realizability model.

\smallskip\noindent
It follows that the axioms of \ ZF are also realized by proof-like $\Cc$-terms in
every realizability model.

\smallskip\noindent
We write the axioms of \ZFe\ with the three binary relation symbols $\neps,\notin,\subseteq$. Of course, \
$x\eps y$ and $x\in y$ \ are the formulas \ $x\neps y\to\bot$ \ and \ $x\notin y\to\bot$.\\
The notation $x\simeq y\to F$ means $x\subseteq y,y\subseteq x\to F$. Thus $x\simeq y$,
which represents the usual (extensional) equality of sets, is the pair of formulas \ $\{x\subseteq y,y\subseteq x\}$.\\
We use the notations \ $(\pt x\eps a)F(x)$ \ for \ $\pt x(\neg F(x)\to x\neps a)$ \ and\\
$(\ex x\eps a)\vec{F}(x)$ for \ $\neg\pt x(\vec{F}(x)\to x\neps a)$.\\
For instance,\ $(\ex x\eps y)\,t\simeq u$ \ is the formula \
$\neg\pt x(t\subseteq u,u\subseteq t\to x\neps y)$.

\smallskip\noindent
The axioms of \ZFe\ are the following~:

\smallskip
0. Extensionality axioms.

$\pt x\pt y[x\in y\dbfl(\ex z\eps y)x\simeq z]$~;
$\pt x\pt y[x\subseteq y\dbfl(\pt z\eps x)z\in y]$.

\smallskip
1. Foundation scheme.

$\pt x_1\ldots\pt x_n\pt a(\pt x((\pt y\eps x)F[y,x_1,\ldots,x_n]\to F[x,x_1,\ldots,x_n])\to
F[a,x_1,\ldots,x_n])$

for every formula $F[x,x_1,\ldots,x_n]$.

\smallskip\noindent
The intuitive meaning of axioms 0 and 1 is that $\varepsilon$ is a well
founded relation, and that the relation $\in$ is obtained by ``~collapsing~''
$\varepsilon$ into an extensional binary relation.
\smallskip\\
The following axioms essentially express that the relation $\varepsilon$
satisfies the axioms of Zermelo-Fraenkel {\em except extensionality}.

\smallskip
2. Comprehension scheme.

$\pt x_1\ldots\pt x_n\pt a\ex b\pt x(x\eps b\dbfl(x\eps a\land F[x,x_1,\ldots,x_n]))$

for every formula $F[x,x_1,\ldots,x_n]$.

\smallskip
3. Pairing axiom.

$\pt a\pt b\ex x\{a\eps x, b\eps x\}$.

\smallskip
4. Union axiom.

$\pt a\ex b(\pt x\eps a)(\pt y\eps x)\,y\eps b$.

\smallskip
5. Power set axiom.

$\pt a\ex b\pt x(\ex y\eps b)
\pt z(z\eps y\dbfl(z\eps a\land z\eps x))$.

\smallskip
6. Collection scheme.

$\pt x_1\ldots\pt x_n\pt a\ex b(\pt x\eps a)(\ex y\,F[x,y,x_1,\ldots,x_n]\to
(\ex y\eps b)F[x,y,x_1,\ldots,x_n])$

for every formula $F[x,y,x_1,\ldots,x_n]$.

\smallskip
7. Infinity scheme.

$\pt x_1\ldots\pt x_n\pt a\ex b\{a\eps b,(\pt x\eps b)
(\ex y\,F[x,y,x_1,\ldots,x_n]\to(\ex y\eps b)F[x,y,x_1,\ldots,x_n])\}$

for every formula $F[x,y,x_1,\ldots,x_n]$.

\medskip\noindent
The usual Zermelo-Fraenkel set theory is obtained from \ZFe\ by identifying the predicate symbols $\neps$ and $\notin$.
Thus, the axioms of ZF are written as follows, with the predicate symbols
$\notin,\subseteq$ (recall that $x\simeq y$ is the conjunction of $x\subseteq y$
and $y\subseteq x)$~:

\smallskip
0. Equality and extensionality axioms.

$\pt x\pt y[x\in y\dbfl(\ex z\in y)x\simeq z]$~;
$\pt x\pt y[x\subseteq y\dbfl(\pt z\in x)z\in y]$.

\smallskip
1. Foundation scheme.

$\pt x_1\ldots\pt x_n\pt a(\pt x((\pt y\in x)F[y,x_1,\ldots,x_n]\to F[x,x_1,\ldots,x_n])\to
F[a,x_1,\ldots,x_n])$

for every formula $F[x,x_1,\ldots,x_n]$ written
with the only relation symbols $\notin,\subseteq$.

\smallskip
2. Comprehension scheme.

$\pt a\ex b\pt x(x\in b\dbfl(x\in a\land F[x,x_1,\ldots,x_n]))$

for every formula $F[x,x_1,\ldots,x_n]$ written
with the only relation symbols $\notin,\subseteq$.

\smallskip
3. Pairing axiom.

$\pt a\pt b\ex x\{a\in x, b\in x\}$.

\smallskip
4. Union axiom.

$\pt a\ex b(\pt x\in a)(\pt y\in x)\,y\in b$.

\smallskip
5. Power set axiom.

$\pt a\ex b\pt x(\ex y\in b)
\pt z(z\in y\dbfl(z\in a\land z\in x))$.

\smallskip
6. Collection scheme.

$\pt x_1\ldots\pt x_n\pt a\ex b(\pt x\in a)(\ex y\,F[x,y,x_1,\ldots,x_n]\to
(\ex y\in b)F[x,y,x_1,\ldots,x_n])$

for every formula $F[x,y,x_1,\ldots,x_n]$ written with the only relation symbols
$\notin,\subseteq$.

\smallskip
7. Infinity scheme.

$\pt x_1\ldots\pt x_n\pt a\ex b\{a\in b,(\pt x\in b)(\ex y\,F[x,y,x_1,\ldots,x_n]\to
(\ex y\in b)F[x,y,x_1,\ldots,x_n])\}$

for every formula $F[x,y,x_1,\ldots,x_n]$ written with the only relation symbols
$\notin,\subseteq$.

\smallskip\noindent
{\small{\bfseries Remark.} The usual statement of the axiom of infinity
is the particular case of this scheme, where $a$ is $\vide$, and $F(x,y)$ is the
formula $y\simeq x\cup\{x\}$.}

\medskip\noindent
Let us show that \ \ZFe\ is a conservative extension of ZF\/.
First, it is clear that, if \ \ZFe\ $\vdash F$, where $F$ is a
formula of ZF (i.e. written only with $\notin$ and $\subseteq$), then \
ZF $\vdash F$~; indeed, it is sufficient to replace $\neps$ with $\notin$
in any proof of \ \ZFe\ $\vdash F$.

\smallskip\noindent
Conversely\/, we must show that each axiom of \ ZF is a consequence of \
\ZFe.

\begin{theorem}\label{a_subset_a}\ \\
\phantom ii)~~\ZFe\ $\vdash\pt a(a\subseteq a)$ (and thus $a\simeq a$).\\
ii)~~\ZFe\ $\vdash\pt a\pt x(x\eps a\to x\in a)$.
\end{theorem}
\proof\hfill\\
\phantom ii)~~Using the foundation axiom, we assume $\pt x(x\eps a\to x\subseteq x)$, and we must show
$a\subseteq a$~; therefore, we add the hypothesis $x\eps a$. It follows that
$x\subseteq x$, then $x\simeq x$, and therefore~:\\
$\ex y\{x\simeq y, y\eps a\}$, that is to say $x\in a$. Thus, we have
$\pt x(x\eps a\to x\in a)$, and therefore $a\subseteq a$.\\
ii)~~Just shown.
\qed

\begin{corollary}\label{b_subset_c}
\ZFe\ $\vdash\pt x(x\in a\to x\in b)\to a\subseteq b$.
\end{corollary}

\proof
We must show $x\eps a\to x\in b$, which follows from $x\in a\to x\in b$ and
$x\eps a\to x\in a$ (theorem~\ref{a_subset_a}(ii)).
\qed

\begin{lemma}\label{a_subset_b_subset_c}
\ZFe\ $\vdash a\subseteq b,\,\pt x(x\in b\to x\in c)\to a\subseteq c$.
\end{lemma}\proof
We must show $x\eps a\to x\in c$, which follows from $x\eps a\to x\in b$ and
$x\in b\to x\in c$.
\qed

\begin{theorem}\label{z_in_a_subset_y}
\ZFe\ $\vdash \pt y\pt z(y\simeq a,\,a\in z\to y\in z)$~; \
\ZFe\ $\vdash\pt y\pt z(a\subseteq y,\,z\in a\to z\in y)$.
\end{theorem}\proof
Call $F(a)$, $F'(a)$ these two formulas. We show $F(a)$ by foundation~:\\
thus, we suppose $(\pt x\eps a)F(x)$ and we first show $F'(a)$~: by hypothesis, we have $a\subseteq y$,
$z\in a$~; thus, there exists $a'$ such that $z\simeq a'$ and $a'\eps a$, and thus
$F(a')$. From $a'\eps a$ and $a\subseteq y$, we deduce $a'\in y$. From $z\simeq a'$
and $a'\in y$, we deduce \ $z\in y$ \ by \ $F(a')$.\\
Then, we show $F(a)$~: by hypothesis, we have $y\simeq a$, $a\in z$, thus $a\simeq y'$ and $y'\eps z$
for some $y'$. In order to show \ $y\in z$, it is sufficient to show $y\simeq y'$.\\
Now, we have $y\simeq a$, $a\simeq y'$, and thus $y'\subseteq a$, $a\subseteq y$. From
$F'(a)$, we get $\pt z(z\in a\to z\in y)$~; from $y'\subseteq a$, we deduce
$y'\subseteq y$ by lemma~\ref{a_subset_b_subset_c}.\\
We have also $y\subseteq a$, $a\subseteq y'$. From $F'(a)$, we get $\pt z(z\in a\to z\in y')$~;
from $y\subseteq a$, we deduce $y\subseteq y'$ by lemma~\ref{a_subset_b_subset_c}.
\qed

\smallskip\noindent
With corollary~\ref{b_subset_c}, we obtain~:
\begin{corollary}
\ZFe\ $\vdash b\subseteq c\dbfl\pt x(x\in b\to x\in c)$.\qed
\end{corollary}\noindent
It is now easy to deduce the equality and extensionality axioms of ZF~:

\smallskip\noindent
$\pt x(x\simeq x)$~; $\pt x\pt y(x\simeq y\to y\simeq x)$~;
$\pt x\pt y\pt z(x\simeq y,y\simeq z\to x\simeq z)$~;\\
$\pt x\pt x'\pt y\pt y'(x\simeq x',y\simeq y',x\notin y\to x'\notin y')$~;
$\pt x\pt y(\pt z(z\notin x\dbfl z\notin y)\to x\simeq y)$~;\\
$\pt x\pt y(x\subseteq y\dbfl\pt z(z\notin y\to z\notin x))$.

\smallskip\noindent
{\small{\bfseries Remark.} This shows that $\simeq$ is an equivalence relation which is compatible
with the relations $\in$ and $\subseteq$~; but, in general, it is
\emph{not compatible with \ $\varepsilon$}. It is the equality relation for ZF~;
it will be called \emph{extensional equivalence}.}

\smallskip\noindent
{\bfseries Notation.} The formula \ $\pt z(z\neps y\to z\neps x)$ \ will be written \ $x\subset y$. 
The ordered pair of formulas \ $x\subset y,y\subset  x$ \ will be written \ $x\sim y$.\\
By theorem~\ref{a_subset_a}, we get \
\ZFe\ $\vdash\pt x\pt y(x\subset y\to x\subseteq y)$. Thus \ $\subset$ \ will be called
\emph{strong inclusion}, and \ $\sim$ \ will be called
\emph{strong extensional equivalence}.

\smallskip\noindent\nopagebreak
$\bullet$~~Foundation scheme.
\smallskip\\
Let $F[x]$ be written with only $\notin,\subseteq$ and let $G[x]$ be the formula
$\pt y(y\simeq x\to F[y])$. Clearly\,, $\pt x\,G[x]$ is
equivalent to $\pt x\,F[x]$. Therefore, from axiom scheme~1 of \ZFe,
it is sufficient to show~: \
$\pt b(\pt x(x\in b\to F[x])\to F[b])\to(\pt x(x\eps a\to G[x])\to G[a])$,
i.e.~:\\
$\pt b(\pt x(x\in b\to F[x])\to F[b]),\pt x\pt y(x\eps a,y\simeq x\to F[y]),
a\simeq b\to F[b]$.\\
Therefore, it is sufficient to prove~: \
$\pt x\pt y(x\eps a,y\simeq x\to F[y]),a\simeq b\to\pt x(x\in b\to F[x])$.\\
From $x\in b,a\simeq b$, we deduce $x\in a$ and therefore (by axiom~0), $x'\eps a$ for some $x'\simeq x$. Finally, we get $F[x]$ from $\pt x\pt y(x\eps a,y\simeq x\to F[y])$.

\smallskip\noindent\label{preuve_axiomes}
$\bullet$~~Comprehension scheme~: $\pt a\ex b\pt x
(x\in b\dbfl(x\in a\land F[x]))$\\
for every formula $F[x,x_1,\ldots,x_n]$ written with $\notin,\subseteq$.\\
From the axiom scheme~2 of \ZFe, we get
$\pt x(x\eps b\dbfl(x\eps a\land F[x]))$. If $x\in b$, then $x\simeq x'$,
$x'\eps b$ for some $x'$. Thus $x'\eps a$ and $F[x']$. From $x\simeq x'$ and
$x'\eps a$, we deduce $x\in a$. Since $\subseteq$ and $\in$ are compatible with $\simeq$,
it is the same for $F$~; thus, we obtain $F[x]$.\\
Conversely\/, if we have $F[x]$ and $x\in a$, we have $x\simeq x'$ and $x'\eps a$
for some $x'$. Since $F$ is compatible with $\simeq$, we get $F[x']$, thus
$x'\eps b$ and $x\in b$.

\smallskip\noindent
$\bullet$~~Pairing axiom~:
$\pt x\pt y\ex z\{x\in z,y\in z\}$.

\smallskip\noindent
Trivial consequence of axiom~3 of \ZFe, and theorem~\ref{a_subset_a}(ii).

\smallskip\noindent
$\bullet$~~Union axiom~: $\pt a\ex b\pt x\pt y(x\in a,y\in x\to y\in b)$.
\smallskip\noindent\\
From $x\in a$ we have $x\simeq x'$ and $x'\eps a$ for some $x'$~; we have $y\in x$,
therefore $y\in x'$, thus $y\simeq y'$ and $y'\eps x'$ for some $y'$. From axiom~4 of \ZFe,
$x'\eps a$ and $y'\eps x'$, we get $y'\eps b$~; therefore $y\in b$, by $y\simeq y'$.

\smallskip\noindent
$\bullet$~~Power set axiom~: $\pt a\ex b\pt x\ex y
\{y\in b,\pt z(z\in y\dbfl(z\in a\land z\in x))\}$

\smallskip\noindent
Given $a$, we obtain $b$ by axiom~5 of \ZFe~; given $x$, we define $x'$ by the
condition~:\\
$\pt z(z\eps x'\dbfl(z\eps a\land z\in x))$ (comprehension scheme of
\ZFe). By definition of $b$, there exists $y\eps b$ such that
$\pt z(z\eps y\dbfl z\eps a\land z\eps x')$, and therefore
$\pt z(z\eps y\dbfl z\eps a\land z\in x)$.\\
It follows easily that $\pt z(z\in y\dbfl z\in a\land z\in x)$.

\smallskip\noindent
$\bullet$~~Collection scheme~: \
$\pt a\ex b(\pt x\in a)(\ex y\,F[x,y]\to(\ex y\in b)F[x,y])$\\
for every formula $F[x,y,x_1,\ldots,x_n]$ written with the only relation symbols
$\notin,\subseteq$.

\smallskip\noindent
From $x\in a$ and $\ex y\,F[x,y]$, we get $x\simeq x'$, $x'\eps a$ for some
$x'$, and thus $\ex y\,F[x',y]$ since $F$ is compatible with~$\simeq$. From
axiom scheme~6 of \ZFe, we get $(\ex y\eps b)F[x',y]$, and therefore
$(\ex y\in b)F[x,y]$, by theorem~\ref{a_subset_a}(ii), again because $F$ is compatible with~$\simeq$.

\smallskip\noindent
$\bullet$~~Infinity scheme~: \
$\pt a\ex b\{a\in b,(\pt x\in b)(\ex y\,F[x,y]\to(\ex y\in b)F[x,y])\}$\\
for every formula $F[x,y,x_1,\ldots,x_n]$ written with the only relation symbols
$\notin,\subseteq$.
\smallskip\\
Same proof.

\section{Realizability models of \ZFe}
As usual in relative consistency proofs, we start with a model ${\mathcal M}$ of ZFC, called \emph{the ground model} \
or \emph{the standard model}. In particular, the integers of ${\mathcal M}$ are called
\emph{the standard integers}.\\
The elements of ${\mathcal M}$ will be called \emph{individuals}.

\smallskip\noindent
In the sequel, the model ${\mathcal M}$ will be our universe, which means that every
notion we consider is defined in ${\mathcal M}$. In particular, the realizability algebra
$(\Lbd,\Pi,\bbot)$ is an individual of ${\mathcal M}$.

\smallskip\noindent
We define a \emph{realizability model}~${\mathcal N}$, with the same set of individuals as
${\mathcal M}$. But ${\mathcal N}$ is not a model in the usual sense, because its truth values
are subsets of $\Pi$ instead of being $0$ or~$1$. Therefore, although ${\mathcal M}$ and
${\mathcal N}$ have the same domain (the quantifier $\pt x$ describes the same domain for
both), the model ${\mathcal N}$ may (and will, in all non trivial cases) have much more
indivi\-duals than~${\mathcal M}$, because it has individuals which are \emph{not named}.
In particular, it will have \emph{non standard integers}.

\smallskip\noindent
{\small{\bfseries Remark.} This is a great difference between \emph{realizability} and \emph{forcing}
models of ZF\/. In a forcing model, each individual is named in the ground model~; it follows that
integers, and even ordinals, are not changed.}

\smallskip\noindent
For each closed formula $F$ with parameters in ${\mathcal M}$, we define two truth values~:\\
$\|F\|\subseteq\Pi$ and $|F|\subseteq\Lbd$.\\
$|F|$ is defined immediately from $\|F\|$ \ as follows~:\\
\centerline{$\xi\in|F|$ \ $\Dbfl$ \ $(\pt\pi\in\|F\|)\,\xi\star\pi\in\bbot$.}

\smallskip\noindent
{\bfseries Notation.} \ We shall write \ $\xi\force F$ \ (read \emph{``~$\xi$ realizes $F$~''}) for \ $\xi\in|F|$.

\smallskip\noindent
$\|F\|$ is now defined by recurrence on the length of $F$~:

\smallskip\noindent
$\bullet$~~$F$ is atomic~;\\
then $F$ has one of the forms \ $\top,\,\bot,\,a\neps b,\,a\subseteq b,\,a\notin b$ where \ $a,b$ are parameters in ${\mathcal M}$.\\
We set~:

\smallskip\noindent
$\|\top\|=\vide$~; \ $\|\bot\|=\Pi$~; \ $\|a\neps b\|=\{\pi\in\Pi;\;(a,\pi)\in b\}$.

\smallskip\noindent
$\|a\subseteq b\|,\|a\notin b\|$ are defined simultaneously by induction on
$($rk$(a)\cup\,$rk$(b),$rk$(a)\cap\,$rk$(b))$\\
(rk$(a)$ being the rank of $a$).

\smallskip\noindent
$\dsp\|a\subseteq b\|=\bigcup_c\{\xi\ps\pi;\;\xi\in\Lbd,\;\pi\in\Pi,\;(c,\pi)\in a,\;
\xi\force c\notin b\}$~;

\smallskip\noindent
$\dsp\|a\notin b\|=\bigcup_c\{\xi\ps\xi'\ps\pi;\;\xi,\xi'\in\Lbd,\;\pi\in\Pi,\;(c,\pi)\in b,\;
\xi\force a\subseteq c,\;\xi'\force c\subseteq a\}$.

\smallskip\noindent
$\bullet$~~$F\equiv A\to B$~; then \
$\|F\|=\{\xi\ps\pi~;\;\xi\force A,\;\pi\in\|B\|\}$.

\smallskip\noindent
$\bullet$~~$F\equiv\pt x\,A$~: then \ $\dsp\|F\|=\bigcup_a\|A[a/x]\|$.

\smallskip\noindent
The following theorem is an essential tool~:

\begin{theorem}[Adequacy lemma]\label{adequat}\ \\
Let $A_1,\ldots,A_n,A$ be closed formulas of \ \ZFe, and suppose that \
$x_1:A_1,\ldots,x_n:A_n\vdash t:A$.\\
If $\xi_1\force A_1,\ldots,\xi_n\force A_n$ then $t[\xi_1/x_1,\ldots,\xi_n/x_n]\force A$.\\
In particular, if $\vdash t:A$, then $t\force A$.
\end{theorem}\noindent
We need to prove a (seemingly) more general result, that we state as a lemma~:

\begin{lemma}\label{adequat_lem}
Let \ $A_1[\vec{z}],\ldots,A_n[\vec{z}],A[\vec{z}]$ be formulas of \ \ZFe, with
$\vec{z}=(z_1,\ldots,z_k)$ as free variables, and suppose that \
$x_1:A_1[\vec{z}],\ldots,x_n:A_n[\vec{z}]\vdash t:A[\vec{z}]$.\\
If \ $\xi_1\force A_1[\vec{a}],\ldots,\xi_n\force A_n[\vec{a}]$ \ for some parameters
(i.e. individuals in~${\mathcal M}$)\\
${\vec{a}=(a_1,\ldots,a_k)}$, then \ $t[\xi_1/x_1,\ldots,\xi_n/x_n]\force A[\vec{a}]$.
\end{lemma}\noindent
\proof By recurrence on the length of the derivation of \
$x_1:A_1[\vec{z}],\ldots,x_n:A_n[\vec{z}]\vdash t:A[\vec{z}]$.\\
We consider the last used rule.

\smallskip\noindent
1. $x_1:A_1[\vec{z}],\ldots,x_n:A_n[\vec{z}]\vdash x_i:A_i[\vec{z}]$. This case
is trivial.

\smallskip\noindent
2. We have the hypotheses~:\\
$x_1:A_1[\vec{z}],\ldots,x_n:A_n[\vec{z}]\vdash u:B[\vec{z}]\to A[\vec{z}]$~~; \
$x_1:A_1[\vec{z}],\ldots,x_n:A_n[\vec{z}]\vdash v:B[\vec{z}]$~~; \ $t=uv$.\\
By the induction hypothesis, we have \
$u[\vec{\xi}/\vec{x}]\force B[\vec{a}/\vec{z}]\to A[\vec{a}/\vec{z}]$ \ and \
$v[\vec{\xi}/\vec{x}]\force B[\vec{a}/\vec{z}]$.\\
Therefore \ $(uv)[\vec{\xi}/\vec{x}]\force A[\vec{a}/\vec{z}]$ \
which is the desired result.

\smallskip\noindent
3. We have the hypotheses~:\\
$x_1:A_1[\vec{z}],\ldots,x_n:A_n[\vec{z}],y:B[\vec{z}]\vdash u:C[\vec{z}]$~~; \
$A[\vec{z}]\equiv B[\vec{z}]\to C[\vec{z}]$~~; \ $t=\lbd y\,u$.\\
We want to show that \
$(\lbd y\,u)[\vec{\xi}/\vec{x}]\force B[\vec{a}/\vec{z}]\to C[\vec{a}/\vec{z}]$.
Thus, let~:\\
$\eta\force B[\vec{a}/\vec{z}]$ \ and \ $\pi\in\|C[\vec{a}/\vec{z}]\|$. We must show~:\\
$(\lbd y\,u)[\vec{\xi}/\vec{x}]\star\eta\ps\pi\in\bbot$ \ or else \
$u[\vec{\xi}/\vec{x},\eta/y]\star\pi\in\bbot$.\\
Now, by the induction hypothesis, we have \
$u[\vec{\xi}/\vec{x},\eta/y]\force C[\vec{a}/\vec{z}]$,\\
which gives the result.

\smallskip\noindent
4. We have the hypotheses~:\\
$x_1:A_1[\vec{z}],\ldots,x_n:A_n[\vec{z}]\vdash t:B[\vec{z}]$~~; \
$A[\vec{z}]\equiv\pt z_1B[\vec{z}]$~~; \
$\xi_i\force A_i[a_1/z_1,a_2/z_2,\ldots,a_k/z_k]$~~;\\
the variable $z_1$ is not free in $A_1[\vec{z}],\ldots,A_n[\vec{z}]$.\\
We have to show that \ $t[\vec{\xi}/\vec{x}]\force\pt z_1B[\vec{a}/\vec{z}]$ \ i.e. \
$t[\vec{\xi}/\vec{x}]\force\pt z_1B[a_2/z_2,\ldots,a_k/z_k]$. Thus, we take an arbitrary
set $b$ in ${\mathcal M}$ and we show \
$t[\vec{\xi}/\vec{x}]\force B[b/z_1,a_2/z_2,\ldots,a_k/z_k]$.\\
By the induction hypothesis, it is sufficient to show that \
$\xi_i\force A_i[b/z_1,a_2/z_2,\ldots,a_k/z_k]$.\\
But this follows from the hypothesis on $\xi_i$, because $z_1$ is not free in the
formulas $A_i$.

\smallskip\noindent
5. We have the hypotheses~:\\
$x_1:A_1[\vec{z}],\ldots,x_n:A_n[\vec{z}]\vdash t:\pt y\,B[y,\vec{z}]$~~; \
$A[\vec{z}]\equiv B[\tau[\vec{z}]/y,\vec{z}]$~~; \ $\xi_i\force A_i[\vec{a}]$.\\
By the induction hypothesis, we have
$t[\vec{\xi}/\vec{x}]\force\pt yB[y,\vec{a}/\vec{z}]$~; therefore
$t[\vec{\xi}/\vec{x}]\force B[b/y,\vec{a}/\vec{z}]$ \ for every parameter $b$.
We get the desired result by taking \ $b=\tau[\vec{a}]$.

\smallskip\noindent
6. The result follows from the following~:

\begin{theorem}
For every formulas $A,B$, we have \ $\Ccc\force ((A\to B)\to A)\to A$.
\end{theorem}

\proof Let \ $\xi\force(A\to B)\to A$ \ and \ $\pi\in\|A\|$. Then \
$\Ccc\star\xi\ps\pi\succ\xi\star\kk_\pi\ps\pi$ \ which is in $\bbot$, because \
$\kk_\pi\force A\to B$ \ by lemma~\ref{k_pi}.
\qed

\begin{lemma}\label{k_pi}
If \ $\pi\in\|A\|$, then \ $\kk_\pi\force A\to B$.
\end{lemma}

\proof Indeed, let \ $\xi\force A$~; then \ $\kk_\pi\star\xi\ps\pi'\succ\xi\star\pi\in\bbot$ \
for every stack \ $\pi'\in\|B\|$.
\qed

\smallskip\noindent
7. We have the hypothesis \ \ $x_1:A_1[\vec{z}],\ldots,x_n:A_n[\vec{z}]\vdash t:\bot$.\\
By the induction hypothesis, we have \ $t[\vec{\xi}/\vec{x}]\force\bot$. Since \ $\|\bot\|=\Pi$,
we have \ $t[\vec{\xi}/\vec{x}]\,\star\,\pi\in\bbot$ for every $\pi\in\|A[\vec{a}/\vec{z}]\|$, and
therefore \ $t[\vec{\xi}/\vec{x}]\force A[\vec{a}/\vec{z}]$ \ which is the desired result.

\smallskip\noindent
This completes the proof of lemma~\ref{adequat_lem} and theorem~\ref{adequat}.
\qed

\subsubsection*{Realized formulas and coherent models}
In the ground model ${\mathcal M}$, we interpret the formulas of the \emph{language of \ ZF}~:
this language consists of $\notin,\subseteq$~; we add some function symbols, but these functions are always
defined, in~${\mathcal M}$, by some formulas written with $\notin,\subseteq$. We suppose that this ground model
satisfies~ZFC.\\
The value, in ${\mathcal M}$, of a closed formula $F$ of the language of ZF, with parameters in ${\mathcal M}$, is of course $1$ or $0$. In the first case, we say that ${\mathcal M}$ \emph{satisfies} $F$, and we write \
${\mathcal M}\models F$.

\smallskip\noindent
In the realizability model ${\mathcal N}$, we interpret the formulas of the \emph{language of \ \ZFe},
which consists of $\neps,\notin,\subseteq$ and the same function symbols as in the language of ZF\/.
The domain of ${\mathcal N}$ and the interpretation of the function symbols are the same as for the
model ${\mathcal M}$.\\
The value, in ${\mathcal N}$, of a closed formula $F$ of \ZFe \ with parameters (in ${\mathcal M}$ or in
${\mathcal N}$, which is the same thing) is an element of ${\mathcal P}(\Pi)$ which is denoted as $\|F\|$,
the definition of which has been given above.\\
Thus, we can no longer say that ${\mathcal N}$ satisfies (or not) a given closed formula $F$. But we shall say
that ${\mathcal N}$ \emph{realizes}~$F$ (and we shall write \ ${\mathcal N}\force F$), if there exists a proof-like term \ $\theta$ such that \ $\theta\force F$. We say that two closed formulas $F\/,G$ are
\emph{interchangeable} if ${\mathcal N}\force F\dbfl G$.\\
Notice that, if $\|F\|=\|G\|$, then $F\/,G$ are interchangeable (indeed \ $I\force F\to G$), but the
converse is far from being true.

\smallskip\noindent
The model ${\mathcal N}$ allows us to make relative consistency proofs, since it is clear, from the adequacy
lemma (theorem~\ref{adequat}), that the class of formulas which are realized in ${\mathcal N}$ is closed by
deduction in classical logic. Nevertheless, we must check that the realizability model ${\mathcal N}$ is
\emph{coherent}, i.e. that it does not realize the formula~$\bot$. We can express this condition in the
following form~:

\smallskip\noindent\centerline{\em For every proof-like term $\theta$, there exists a stack $\pi\in\Pi$
such that $\theta\star\pi\notin\bbot$.}

\smallskip\noindent
When the model ${\mathcal N}$ is coherent, it is not \emph{complete}, except in trivial cases.
This means that there exist closed formulas $F$ of \ZFe \ such that
${\mathcal N}\nforce F$ \ and \ ${\mathcal N}\nforce\neg F$.

\subsection*{The axioms of \ZFe \ are realized in ${\mathcal N}$}\ 

\bigskip\noindent
$\bullet$~~Extensionality axioms.

\smallskip\noindent
We have \ $\dsp\|\pt z(z\notin b\to z\neps a)\|
=\bigcup_c\{\xi\ps\pi;\;\xi\force c\notin b,\,\pi\in\|c\neps a\|\}$\\
by definition of the value of $\|\pt z(z\notin b\to z\neps a)\|$~;

\smallskip\noindent
and \ $\dsp\|a\subseteq b\|=\bigcup_c\{\xi\ps\pi;\;(c,\pi)\in a,\;\xi\force c\notin b\}$ \
by definition of \ $\|a\subseteq b\|$.\\
Therefore, we have \ $\|a\subseteq b\|=\|\pt z(z\notin b\to z\neps a)\|$, so that~:\\
$I\force\pt x\pt y(x\subseteq y\to\pt z(z\notin y\to z\neps x))$ \ and \
$I\force\pt x\pt y(\pt z(z\notin y\to z\neps x)\to x\subseteq y)$.

\smallskip\noindent
In the same way, we have~:

\smallskip\noindent
$\dsp\|\pt z(a\subseteq z,z\subseteq a\to z\neps b)\|=\bigcup_c\{\xi\ps\xi'\ps\pi;\;
\xi\force a\subseteq c,\;\xi'\force c\subseteq a;\;\pi\in\|c\neps b\|\}$\\
by definition of the value of \ $\|\pt z(a\subseteq z,z\subseteq a\to z\neps b)\|$~;

\smallskip\noindent
and \ $\dsp\|a\notin b\|=\bigcup_c
\{\xi\ps\xi'\ps\pi;\;(c,\pi)\in b,\;\xi\force a\subseteq c,\;\xi'\force c\subseteq a\}\}$ \
by definition of \ $\|a\notin b\|$.\\
Therefore, we have \ $\|a\notin b\|=\|\pt z(a\subseteq z,z\subseteq a\to z\neps b)\|$,
so that~:\\
$I\force\pt x\pt y(x\notin y\to\pt z(x\subseteq z,z\subseteq x\to z\neps y))$~;\\
$I\force\pt x\pt y(\pt z(x\subseteq z,z\subseteq x\to z\neps y)\to x\notin y)$.

\medskip\noindent
{\bfseries Notation.} We shall write $\vec{\xi}$ for a finite sequence $(\xi_1,\ldots,\xi_n)$
of terms. Therefore, we shall write $\vec{\xi}\force\vec{A}$ for $\xi_i\force A_i$
($i=1,\ldots,n$).\\
In particular, the notation $\vec{\xi}\force a\simeq b$ means
$\xi_1\force a\subseteq b,\,\xi_2\force b\subseteq a$~;\\
the notation $\vec{\xi}\force A\dbfl B$ means $\xi_1\force A\to B$, $\xi_2\force B\to A$.

\bigskip\noindent
$\bullet$~~Foundation scheme.

\begin{theorem}\label{fondation}
For every finite sequence $\vec{F}[x,x_1,\ldots,x_n]$ of formulas, we have~:\\
$\Y\force\pt x(\pt y(\vec{F}[y]\to y\neps x),\vec{F}[x]\to\bot)\to
\pt x(\vec{F}[x]\to\bot)$\\
with \ $\Y=AA$ \ and \ $A=\lbd a\lbd f(f)(a)af$ (Turing fixed point combinator).
\end{theorem}

\proof Let $\xi\force\pt x(\pt y(\vec{F}[y]\to y\neps x),\vec{F}[x]\to\bot)$. We show, by induction
on the rank of $a$, that~:\\
$\Y\star\xi\ps\vec{\eta}\ps\pi\in\bbot$, for every \ $\pi\in\Pi$ and \ $\vec{\eta}\force\vec{F}[a]$.\\
Since \ $\Y\star\xi\ps\vec{\eta}\ps\pi\succ\xi\star\Y\xi\ps\vec{\eta}\ps\pi$, it suffices to show \ $\xi\star\Y\xi\ps\vec{\eta}\ps\pi\in\bbot$.\\
Now, \ $\xi\force\pt y(\vec{F}[y]\to y\neps a),\vec{F}[a]\to\bot$, so that it suffices
to show \ $\Y\xi\force\pt y(\vec{F}[y]\to y\neps a)$, in other words \
$\Y\xi\force\vec{F}[b]\to b\neps a$ \ for every $b$. Let $\vec{\zeta}\force\vec{F}[b]$
and $\varpi\in\|b\neps a\|$. Thus, we have $(b,\varpi)\in a$, therefore rk$(b)<$ rk$(a)$ so
that $\Y\star\xi\ps\vec{\zeta}\ps\varpi\in\bbot$ \ by induction hypothesis. It follows
that \ $\Y\xi\star\vec{\zeta}\ps\varpi\in\bbot$, which is the desired result.\qed

\smallskip\noindent
It follows from theorem~\ref{fondation} that the axiom scheme~1 of \ZFe\ (foundation) is realized.

\bigskip\noindent
$\bullet$~~Comprehension scheme.\\
Let $a$ be a set, and $F[x]$ a formula with parameters. We put~:\\
$b=\{(x,\xi\ps\pi);$ $(x,\pi)\in a$, $\xi\force F[x]\}$~; then, we have trivially \
$\|x\neps b\|=\|F(x)\fl x\neps a\|$.\\
Therefore $I\force\pt x(x\neps b\to(F(x)\to x\neps a))$ \ and \
$I\force\pt x((F(x)\to x\neps a)\to x\neps b)$.

\bigskip\noindent
$\bullet$~~Pairing axiom.\\
We consider two sets $a$ and $b$, and we put $c=\{a,b\}\fois\Pi$. We have \
$\|a\neps c\|=\|b\neps c\|=\|\bot\|$, thus $I\force a\eps c$ and $I\force b\eps c$.\\
{\small{\bfseries Remark.}\\
Except in trivial cases, $c$ has many other elements than $a$ and $b$, which have no name in
${\mathcal M}$.}

\bigskip\noindent
$\bullet$~~Union axiom.\\
Given a set $a$, let $b=$ Cl$(a)$ (the transitive closure of $a$, i.e. the least transitive set
which contains~$a$). We show \ $\|y\neps b\to x\neps a\|\subseteq\|y\neps x\fl x\neps a\|$~:\\
indeed, let $\xi\ps\pi\in\|y\neps b\to x\neps a\|$, i.e. \ $\xi\force y\neps b$ \ and \
$(x,\pi)\in a$.\\
Therefore, \ $x\subseteq$ Cl$(a)$, i.e. \ $x\subseteq b$ \ and thus \
$\|y\neps b\|\supset\|y\neps x\|$.\\
Thus, we have \ $\xi\force y\neps x$, which gives the result.\\
It follows that $I\force\pt x\pt y((y\neps x\fl x\neps a)\to(y\neps b\fl x\neps a))$.

\bigskip\noindent
$\bullet$~~Power set axiom.\\
Given a set $a$, let $b={\mathcal P}($Cl$(a)\fois\Pi)\fois\Pi$. For every set $x$, we put~:\\
$y=\{(z,\xi\ps\pi);$ $\xi\force z\eps x,\,(z,\pi)\in a\}$. We have
$y=\{(z,\xi\ps\pi);$ $\xi\force z\eps x,\,\pi\in\|z\neps a\|\}$, and therefore \
$\|z\neps y\|=\|z\eps x\to z\neps a\|$. Thus~:\\
$I\force\pt z(z\neps y\to(z\eps x\to z\neps a))$ \
and \ $I\force\pt z((z\eps x\to z\neps a)\to z\neps y)$.\\
Now, it is obvious that $y\in{\mathcal P}($Cl$(a)\fois\Pi)$, and therefore
$(y,\pi)\in b$ for every $\pi\in\Pi$.\\
Thus, we have \ $\|y\neps b\|=\Pi=\|\bot\|$. It follows that~:\\
$\lbd f(f)II\force\pt x(\pt y(\pt z(z\neps y\to(z\eps x\to z\neps a)),
\pt z((z\eps x\to z\neps a)\to z\neps y)\to y\neps b)\to\bot)$.

\bigskip\noindent
$\bullet$~~Collection scheme.\\
Given a set $a$, and a formula $F[x,y]$ with parameters, let~:

$b=\bigcup\{\Phi(x,\xi)\fois\,$Cl$(a)$; $x\in\,$Cl$(a),\xi\in\Lbd\}$ \ with

$\Phi(x,\xi)=\{y$ of minimum rank $;\;\xi\force F[x,y]\}$ \ or \ $\Phi(x,\xi)=\vide$ \
if there is no such~$y$.

\noindent
We show that \
$\|\pt y(F[x,y]\fl x\neps a)\|\subseteq\|\pt y(F[x,y]\fl y\neps b)\|$~:\\
Suppose indeed that \ $\xi\ps\pi\in\|\pt y(F[x,y]\fl x\neps a)\|$, i.e. \ $(x,\pi)\in a$ \
and \ $\xi\force F[x,y]$ \ for some~$y$.
By definition of $\Phi(x,\xi)$, there exists $y'\in\Phi(x,\xi)$. Moreover, we have~:\\
$x\in$ Cl$(a),\,\pi\in$ Cl$(a)$, and therefore $(y',\pi)\in b$~; it follows that
$\pi\in\|y'\neps b\|$. But, since $y'\in\Phi(x,\xi)$, we have
$\xi\force F[x,y']$ \ and thus \ $\xi\ps\pi\in\|F[x,y']\to y'\neps b\|$, which gives the result.
We have proved that \
$I\force\pt x(\pt y(F[x,y]\fl y\neps b)\to\pt y(F[x,y]\fl x\neps a))$.

\bigskip\noindent
$\bullet$~~Infinity scheme.\\
Given a set $a$, we define $b$ as the least set such that~:\\
\centerline{$\{a\}\fois\Pi\subseteq b$ \ and \
$\pt x(\pt\pi\in\Pi)(\pt\xi\in\Lbd)((x,\pi)\in b\,\Fl\Phi(x,\xi)\fois\{\pi\}\subseteq b)$}\\
where \ $\Phi(x,\xi)$ \ is defined as above.\\
We have \ $\{a\}\fois\Pi\subseteq b$, thus \ $\|a\neps b\|=\|\bot\|$, and therefore \
$I\force a\eps b$.\\
We now show that \ $\|\pt y(F[x,y]\fl x\neps b)\|\subseteq\|\pt y(F[x,y]\fl y\neps b)\|$~:\\
Suppose indeed that \ $\xi\ps\pi\in\|\pt y(F[x,y]\fl x\neps b)\|$, i.e. \ $(x,\pi)\in b$ \
and \ $\xi\force F[x,y]$ \ for some~$y$.
By definition of $\Phi(x,\xi)$, there exists $y'\in\Phi(x,\xi)$. By definition of $b$, we have \
$(y',\pi)\in b$, \ i.e. \ $\pi\in\|y'\neps b\|$. Now, since $y'\in\Phi(x,\xi)$, we have
$\xi\force F[x,y']$ \ and thus~:\\
$\xi\ps\pi\in\|F[x,y']\to y'\neps b\|$, which gives the result.\\
We have proved that \ $I\force a\eps b$ \ and \
$I\force\pt x(\pt y(F[x,y]\fl y\neps b)\to\pt y(F[x,y]\fl x\neps b))$.

\subsection*{Function symbols and equality}\ \\
According to our needs, we shall add to the language of \ZFe, some \emph{function symbols}
$f,g,\ldots$ \ of any arity. A $k$-ary function symbol $f$ will be interpreted, in the realizability
model ${\mathcal N}$, by a \emph{functional relation}, which is defined
\emph{in the ground model ${\mathcal M}$} by a formula \ $F[x_1,\ldots,x_k,y]$ of~ZF\/. Thus, we assume
that \ ${\mathcal M}\models \pt x_1\ldots\pt x_k\ex!y\,F[x_1,\ldots,x_k,y]$\\
($\ex!y\,F[y]$ \ is the conjunction of \ $\pt y\pt y'(F[y],F[y']\to y=y')$ \ and \
$\ex y\,F[y]$).\\
The axiom schemes of \ZFe, written in the extended language, are still realized in the model ${\mathcal N}$, because the above proofs remain valid.\\
On the other hand, in order to make sure that the axiom schemes of ZF\/, which use a $k$-ary function
symbol~$f$, are still realized, one must check that this symbol is \emph{compatible with $\simeq$}, i.e. that the following formula is realized in ${\mathcal N}$~:\\
$\pt x_1\ldots\pt x_k(x_1\simeq y_1,\ldots,x_k\simeq y_k\to fx_1\ldots x_k\simeq fy_1\ldots y_k)$.

\smallskip\noindent
We now add a new rule to build formulas of \ZFe~:

\smallskip\noindent
If $t,u$ are two terms and \ $F$ \ is a formula of \ZFe, then \ $t=u\mapto F$ \ is a formula of \ZFe.

\smallskip\noindent
The formula \ $t=u\mapto\bot$ \ is denoted \ $t\ne u$.\\
The formula \ $t\ne u\to\bot$, i.e. \ $(t=u\mapto\bot)\to\bot$ \ is denoted \ $t=u$.

\smallskip\noindent
The truth value of these new formulas is defined as follows, assuming that \ $t,u,F$ \ are closed,
with parameters in ${\mathcal N}$~:

\smallskip\noindent
$\|t=u\mapto F\|=\vide$ \ if $t\ne u$~; \ $\|t=u\mapto F\|=\|F\|$ \ if $t=u$.

\smallskip\noindent
It follows that~:\\
$\|t\ne u\|=\vide=\|\top\|$ \ if $t\ne u$~; $\|t\ne u\|=\Pi=\|\bot\|$ \ if $t=u$~;\\
$\|t=u\|=\|\top\to\bot\|$ \ if $t\ne u$~; $\|t=u\|=\|\bot\to\bot\|$ \ if $t=u$.

\smallskip\noindent
Proposition~\ref{mapsto} shows that \ $t=u\mapto F$ \ and \
$t=u\to F$ \ are interchangeable.

\begin{proposition}\label{mapsto}\ \\
\phantom ii)~~$\lbd x(x)I\force(t=u\to F)\to(t=u\mapto F)$~;\\
ii)~~$\lbd x\lbd y(\Ccc)\lbd k(y)(k)x\force(t=u\mapto F),t=u\to F$.
\end{proposition}
\proof\hfill\\
\noindent\phantom ii)~Let $\xi\force t=u\to F$ and $\pi\in\|t=u\mapto F\|$. Thus, we have $t=u$ and $\pi\in\|F\|$.\\
We must show $\lbd x(x)I\star\xi\ps\pi\in\bbot$, that is $\xi\star I\ps\pi\in\bbot$.
This is immediate, by hypothesis on $\xi$, since $I\force t=u$.

\smallskip\noindent
ii)~Let \ $\xi\force t=u\mapto F$, \ $\eta\force t=u$ \ and \ $\pi\in\|F\|$. We must show that~:\\
$\lbd x\lbd y(\Ccc)\lbd k(y)(k)x\star\xi\ps\eta\ps\pi\in\bbot$, soit \ $\eta\star\kk_\pi\xi\ps\pi\in\bbot$.\\
If $t\ne u$, then \ $\eta\force\top\to\bot$, hence the result.\\
If $t=u$, then $\xi\force F$, thus $\xi\star\pi\in\bbot$, therefore $\kk_\pi\xi\force\bot$.\\
But we have $\eta\force\bot\to\bot$, and therefore $\eta\star\kk_\pi\xi\ps\pi\in\bbot$.
\qed

\smallskip\noindent
Proposition~\ref{leibniz} shows that the formulas \ $t=u$ \ and \
$\pt x(u\neps x\to t\neps x)$ \emph{(Leibniz equality)} \ are interchangeable. 

\begin{proposition}\label{leibniz}\ \\
\phantom ii)~~$I\force t=u\mapto\pt x(u\neps x\to t\neps x)$~;\\
ii)~~$I\force\pt x(u\neps x\to t\neps x)\to t=u$.
\end{proposition}
\proof\hfill\\
\noindent\phantom ii)~~It suffices to check that $I\force\pt x(u\neps x\to t\neps x)$ when $t=u$,
which is obvious.\\
ii)~~We must show that \ $I\force\pt x(u\neps x\to t\neps x),t\ne u\to\bot$. Thus
let $\xi\force\pt x(u\neps x\to t\neps x)$, $\eta\force t\ne u$ and $\pi\in\Pi$~;
we must show that \ $\xi\star\eta\ps\pi\in\bbot$.\\
We have \ $\xi\force u\neps a\to t\neps a$ for every $a$~; we take $a=\{t\}\fois\Pi$,
thus $\|t\neps a\|=\Pi$, hence $\pi\in\|t\neps a\|$.\\
If $t=u$, we have $\eta\force\bot$, thus \ $\eta\force u\neps a$, hence the result.\\
If $t\ne u$, we have $\|u\neps a\|=\vide=\|\top\|$, thus \ $\eta\force u\neps a$, hence the result.
\qed

\smallskip\noindent
We now show that the axioms of equality are realized.

\begin{proposition}\label{ax_eg}
$I\force\pt x(x=x)$~; $I\force\pt x\pt y(x=y\mapto y=x)$~;\\
$I\force\pt x\pt y\pt z(x=y\mapto(y=z\mapto x=z))$~;\\
$I\force\pt x\pt y(x=y\mapto(F[x]\to F[y]))$ for every formula $F$ with one free variable, with parameters.
\end{proposition}

\proof
Trivial, by definition of \ $\mapto$.
\qed

\subsubsection*{Conservation of well-foundedness}
Theorem~\ref{bien_fonde} says that every well founded relation in the ground model ${\mathcal M}$, gives a
well foun\-ded relation in the realizability model ${\mathcal N}$.

\begin{theorem}\label{bien_fonde}
Let $f$ be a binary function such that $f(x,y)=1$ is a well founded relation in the ground model
${\mathcal M}$.
Then, for every formula $F[x]$ of \ \ZFe \ with parameters in ${\mathcal M}$~:\\
$\Y\force\pt x(\pt y(f(y,x)=1\mapto F[y])\to F[x])\to\pt x\,F[x]$\\
with \ $\Y=AA$ \ and \ $A=\lbd a\lbd f(f)(a)af$.
\end{theorem}

\proof
Let us fix $a$ and let \ $\xi\force\pt x(\pt y(f(y,x)=1\mapto F[y])\to F[x])$. We show\/,
by induction on $a$, following the well founded relation $f(x,y)=1$, that $\Y\star\xi\ps\pi\in\bbot$
for every $\pi\in\|F[a]\|$.\\
Thus, suppose that \ $\pi\in\|F[a]\|$~; since \ $\Y\star\xi\ps\pi\succ\xi\star\Y\xi\ps\pi$,
we need to show that \ $\xi\star\Y\xi\ps\pi\in\bbot$.
By hypothesis, we have \ $\xi\force\pt y(f(y,a)=1\mapto F[y])\to F[a]$.\\
Thus, it suffices to show that \ $\Y\xi\force f(y,a)=1\mapto F[y]$ for every $y$.\\
This is clear if $f(y,a)\ne1$, by definition of~~$\mapto$.\\
If $f(y,a)=1$, we must show \ $\Y\xi\force F[y]$, \ i.e. \ $\Y\star\xi\ps\rho\in\bbot$ for every \
$\rho\in\|F[y]\|$. But this follows from the induction hypothesis.
\qed

\subsection*{Sets in ${\mathcal M}$ give type-like sets in ${\mathcal N}$}\ \\
\smallskip\noindent
We define a unary function symbol $\gimel$ \ by putting \ $\gimel(a)=a\fois\Pi$ \ for
every individual $a$ (element of the ground model ${\mathcal M}$).

\smallskip\noindent
For each set $E$ of the ground model ${\mathcal M}$, we also introduce the unary function
$1_E$ with values in \ $\{0,1\}$, defined as follows~:\\
$1_E(a)=1$ if $a\in E$~; $1_E(a)=0$ if $a\notin E$.\\
The formula \ $1_E(x)=1\mapto A$ \ will also be denoted as \ $x\eps\gimel E\mapto A$.\\
In particular, \ $a\neps\gimel E$ \ is identical with \ $a\eps\gimel E\mapto\bot$
that is \ $1_E(a)\ne1$.
\smallskip\\
We shall write \ $\pt x^{\gimel E}A[x]$ \ for \ $\pt x(x\eps\gimel E\mapto A[x])$.\\
Proposition~\ref{mapsto} shows that \ $x\eps\gimel E\mapto A$ \ and \
$x\eps\gimel E\to A$ \ are interchangeable.\\
Therefore \ $\pt x^{\gimel E}A[x]$ \ and \ $\pt x(x \eps\gimel E\to A[x])$ \
are also interchangeable. We have~:

\medskip\noindent
\centerline{$\dsp\|\pt x^{\gimel E}A[x]\|=\bigcup_{a\in E}\|A[a/x]\|$ \ and \
$\dsp|\pt x^{\gimel E}A[x]|=\bigcap_{a\in E}|A[a/x]|$.}

\medskip\noindent
As already said, we shall add to the language of \ \ZFe, some function symbols of any arity,
which will be interpreted in the ground model ${\mathcal M}$ by some functional relations.
Then every formula of the form
$\pt\vec{x}(t_1[\vec{x}]=u_1[\vec{x}],\ldots,t_k[\vec{x}]=u_k[\vec{x}]\to t[\vec{x}]=u[\vec{x}])$
which is satisfied in the model ${\mathcal M}$, is \emph{realized} in the model ${\mathcal N}$
($t_1,u_1,\ldots,t_k,u_k,t,u$ are terms of the language).\\
Indeed, we verify immediately that~:\\
$I\force\pt\vec{x}(t_1[\vec{x}]=u_1[\vec{x}]\mapto(\ldots\mapto(t_k[\vec{x}]=u_k[\vec{x}]\mapto
t[\vec{x}]=u[\vec{x}]))\ldots)$.\\
It follows that if, for instance, $t[x_0,x_1]$ sends $E_0\fois E_1$ into $D$ in the model
${\mathcal M}$, then it sends \ $\gimel E_0\fois\gimel E_1$ \ into \ $\gimel D$ in the model ${\mathcal N}$. Indeed, we have then~:\\
${\mathcal M}\models\pt x_0\pt x_1(1_{E_0}(x_0)=1,1_{E_1}(x_1)=1\to 1_D(t[x_0,x_1])=1)$ and therefore,
we have~:\\
$I\force\pt x_0\pt x_1(1_{E_0}(x_0)=1\mapto(1_{E_1}(x_1)=1\mapto 1_D(t[x_0,x_1])=1))$, in other words~:
\smallskip\\
$I\force\pt x_0^{\gimel E_0}\pt x_1^{\gimel E_1}(t[x_0,x_1]\eps\gimel D)$.

\smallskip\noindent
Notice, in particular, that the characteristic function $1_E$, which takes its values
in the set $\mathbf{2}=\{0,1\}$ in the model ${\mathcal M}$,  sends $\gimel E$ into $\gd$ in the realizability model ${\mathcal N}$.

\smallskip\noindent
We shall denote \ $\et,\ou,\non$ \ the (trivial) Boolean algebra operations in
$\{0,1\}$ (they should not be confused with the logical connectives $\land,\lor,\neg$).
In this way, we have defined three function symbols of the language of \ZFe~;
thus, in the realizability model ${\mathcal N}$, they define a
\emph{Boolean algebra structure} on the set $\gd$.

\smallskip\noindent
{\small{\bfseries Remarks.}\\
\phantom ii)~~A set of the form \ $\gimel E$ \ behaves somewhat like a \emph{type}, in the sense of
computer science, because any function of the model ${\mathcal M}$ with domain (resp. range)
$E_1\fois\cdots\fois E_k$ becomes a function of the model ${\mathcal N}$ with domain (resp. range) \
$\gimel E_1\fois\cdots\fois\gimel E_k$.\\
ii)~~The Boolean algebra $\gd$ is, in general, non trivial i.e. it has $\varepsilon$-elements $\ne0,1$.
Notice that they are all empty~: indeed, it is easy to check that \
$I\force\pt x^{\gd}\pt y(x\ne1\to y\neps x)$.}

\subsection*{The set \ $\wt{\NN}$ \ of integers in ${\mathcal N}$}\ \\
We add to the language of \ZFe \ a constant symbol $0$ and a unary function symbol $s$.
Their interpretation in the model ${\mathcal M}$ is as follows~:\\
$0$ \ is \ $\vide$~; \ $s(a)$ \ is \ $\{a\}\fois\Pi$ \ for every set $a$, in other words \
$s(a)=\gimel(\{a\})$.\\
In the realizability model ${\mathcal N}$, $s(a)$ is the singleton of $a$. Indeed, we have trivially~:\\
$\|b\neps s(a)\|=\|b\ne a\|$ (i.e. $\vide$ if $a\ne b$ and $\Pi$ if $a=b$) and it follows that~:\\
$I\force\pt x\pt y(y\neps sx\to x\ne y)$~; $I\force\pt x\pt y(x\ne y\to y\neps sx)$.

\smallskip\noindent
For each $n\in\NN$, the term $s^n0$ will also be written $n$.

\smallskip\noindent
{\small{\bfseries Remark.} In the definition of the set of integers in the realizability
model ${\mathcal N}$, we prefer to use the singleton as the successor function $s$, instead of the usual one \
$x\longmapsto x\cup\{x\}$, which is more complicated to define. It would give~:
$s(a)=\{(a,K\!\ps\pi);\;\pi\in\Pi\}\cup\{(x,\ul{0}\ps\pi);\;(x,\pi)\in a\}$.}

\begin{theorem}\label{symboles_0_s}
The following formulas are realized in ${\mathcal N}$~:\\
\emph{\phantom{ii}i)}~~$\pt x\pt y(sx=sy\mapto x=y)$~;\\
\emph{\phantom iii)}~~$\pt x(sx\not\simeq 0)$~;\\
\emph{iii)}~~$\pt x\pt y(x\simeq y\to sx\simeq sy)$~;\\
\emph{iv)}~~$\pt x\pt y(sx\simeq sy\to x\simeq y)$.
\end{theorem}\noindent
This shows, in particular, that the function $s$ is \emph{compatible with the extensional
equivalence~$\simeq$}.

\proof\hfill\\
\noindent\phantom{ii}i)~~We check that \ $I\force sa=sb\mapto a=b$. We may suppose $sa=sb$, because\\
$\|sa=sb\mapto a=b\|=\vide$ \ if \ $sa\ne sb$. But, in this case, we have $a=b$, by definition of
$sa,sb$.

\smallskip\noindent
\phantom iii)~~We have $\|a\notin0\|=\|\pt x(x\simeq a\to x\neps0)\|=\vide$, since $\|x\neps0\|=\vide$.
Now $\|a\neps sa\|=\Pi$ and therefore we have, for any $\xi\in\Lbd$, \
$\lbd x(x)\xi\force(a\notin\vide\to a\neps sa)\to\bot$~; thus~:\\
$\lbd x(x)\xi\force\pt x(x\notin\vide\to x\neps sa)\to\bot$. But this means exactly that \
$\lbd x(x)\xi\force sa\subseteq0\to\bot$, and therefore
$\lbd x\lbd y(x)\xi\force sa\simeq0\to\bot$.

\smallskip\noindent
iii)~~We show that the formula \ $a\simeq b\to sa\simeq sb$ is realized~; it suffices to realize the
formula \ $a\simeq b\to sa\subseteq sb$. We prove it by means of already realized sentences.\\
We need to prove \ $a\simeq b,x\notin sb\to x\neps sa$. But $x\neps sa$ has the same truth value as
$x\ne a$. Thus, we simply have to prove $a\simeq b\to a\in sb$. But $a\in sb$ follows from $b\eps sb$
and $a\simeq b$.

\smallskip\noindent
iv)~~In the same way, we prove the formula \ $sa\simeq sb\to a\simeq b$ \ and, in fact \
$sa\subseteq sb\to a\simeq b$.\\
The formula \ $sa\subseteq sb$ \ is \ $\pt x(x\notin sb\to x\neps sa)$~; but \ $x\neps sa$ \ is the
same as $x\ne a$. Thus, from $sa\subseteq sb$ we obtain  $a\in sb$, \ i.e. \ $(\ex x\eps sb)\,x\simeq a$.
But \ $x\eps sb$ \ is the same as $x=b$, \ so that we obtain \ $a\simeq b$.
\qed

\smallskip\noindent
The individuals $s^n0$ are obviously distinct, for $n\in\NN$. Therefore, we can define~:

\smallskip
\centerline{$\wt{\NN}=\{(s^n0,\ul{n}\ps\pi);\;n\in\NN,\,\pi\in\Pi\}$}

\noindent
and we have~:\\
$\|a\neps\wt{\NN}\|=\vide$ if $a$ is not of the form $s^n0$, with $n\in\NN$~;\\
$\|s^n0\neps\wt{\NN}\|=\{\ul{n}\ps\pi;\;\pi\in\Pi\}$.\\
The formula \ $x\eps\wt{\NN}$ will also be written \ ent$(x)$.\\
In the sequel, we shall use the restricted quantifier \ $\pt x^{\wt{\NN}}$, which we also write \
$\pt x\inde$, \ with the following meaning~:\\
$\|\pt x\inde F[x]\|=\|\pt x^{\wt{\NN}}F[x]\|=\{\ul{n}\ps\pi;\;n\in\NN,\;\pi\in\|F[s^n0]\|\}$.\\
The restricted existential quantifier \ $\ex x^{\wt{\NN}}$ \ or \ $\ex x\inde$ is defined as~:\\
$\ex x\inde F[x]\equiv\ex x^{\wt{\NN}}\,F[x]\equiv\neg\pt x\inde\neg F[x]$.\\
Proposition~\ref{ptxn} shows that these quantifiers have indeed the intended meaning~: the formulas \
$\pt x\inde\,F[x]$ \ and \ $\pt x(x\eps\wt{\NN}\to F[x])$ are interchangeable.

\begin{proposition}\label{ptxn}\ \\
\emph{\phantom ii)}~~$\lbd x\lbd y\lbd z(y)(x)z\force\pt x\inde\,F[x]\to\pt x(\neg F[x]\to x\neps\wt{\NN})$~;\\
\emph{ii)}~~$\lbd x\lbd y(\Ccc)\lbd k(x)ky\force\pt x(\neg F[x]\to x\neps\wt{\NN})\to\pt x\inde\,F[x]$.
\end{proposition}

\proof\hfill\\
\noindent\phantom ii)~~Let $\xi\force\pt x\inde\,F[x]$, $\eta\force\neg F[a]$ \ and \ $\varpi\in\|a\neps\wt{\NN}\|$.
Thus, we have \
$a=s^n0$ for some $n\in\NN$ (else $\|a\neps\wt{\NN}\|=\vide$) and \ $\varpi=\ul{n}\ps\pi$.
We must show that \ $\eta\star\xi\ul{n}\ps\pi\in\bbot$.\\
Now, by hypothesis on $\xi$, we have \ $\xi\star\ul{n}\ps\rho\in\bbot$
for any \ $\rho\in\|F[s^n0]\|$~; thus \ $\xi\ul{n}\force F[s^n0]$. Since \ $\eta\force\neg F[s^n0]$,
we have \ $\eta\star\xi\ul{n}\ps\pi\in\bbot$, which is the desired result.

\smallskip\noindent
ii)~~Let \ $\xi\force\pt x(\neg F[x]\to x\neps\wt{\NN})$ \ and \ $\ul{n}\ps\pi\in\|\pt x\inde\,F[x]\|$,
with \ $n\in\NN$ \ and \ $\pi\in\|F[s^n0]\|$.\\
We have~: \
$\lbd x\lbd y(\Ccc)\lbd k(x)ky\star\xi\ps\ul{n}\ps\pi\succ\xi\star\kk_\pi\ps\ul{n}\ps\pi$.

\smallskip\noindent
Now, we have \ $\kk_\pi\force\neg F[s^n0]$ \ and \ $\ul{n}\ps\pi\in\|s^n0\neps\wt{\NN}\|$.
Therefore \ $\xi\star\kk_\pi\ps\ul{n}\ps\pi\in\bbot$.
\qed

\begin{theorem}[Recurrence scheme]\label{rec_eps}
For every formula $F[\vec{x},y]$~:\\
\emph{\phantom ii)}~~$I\force\pt\vec{x}\,\pt n^{\wt{\NN}}\,(\pt y(F[\vec{x},sy]\to F[\vec{x},y]),F[\vec{x},n]\to
F[\vec{x},0])$.\\
\emph{ii)}~~$I\force\pt\vec{x}\,\pt n^{\wt{\NN}}\,(\pt y(F[\vec{x},y]\to F[\vec{x},sy]),
F[\vec{x},0]\to F[\vec{x},n])$.
\end{theorem}

\proof\hfill\\
\noindent\phantom ii)~~Let $n\in\NN$, \ $\vec{a}$ \ a sequence of individuals,
$\xi\force\pt y(F[\vec{a},sy]\to F[\vec{a},y])$, $\pi\in\|F[\vec{a},0]\|$.\\
We must show that, for every \ $\alpha\force F[\vec{a},n]$, \ we have
$I\star\ul{n}\ps\xi\ps\alpha\ps\pi\in\bbot$.\\
In fact, we show, by recurrence on $n$, that \ $\ul{n}\star\xi\ps\alpha\ps\pi\in\bbot$.\\
This is immediate if $n=0$. In order to go from $n$ \ to $n+1$, we suppose now \
$\alpha\force F[\vec{a},sn]$~;\\
we have \ $\ul{n+1}\star\xi\ps\alpha\ps\pi\succ\sig\ul{n}\star\xi\ps\alpha\ps\pi\succ
\sigma\star\ul{n}\ps\xi\ps\alpha\ps\pi\succ\ul{n}\star\xi\ps\xi\alpha\ps\pi$.\\
But, by hypothesis on $\xi$, we have \
$\xi\force F[\vec{a},sn]\to F[\vec{a},n]$~; thus \ $\xi\alpha\force F[\vec{a},n]$.\\
Hence the result, by the recurrence hypothesis.

\smallskip\noindent
ii)~~Let $n\in\NN$, \ $\vec{a}$ \ a sequence of individuals,
$\xi\force\pt y(F[\vec{a},y]\to F[\vec{a},sy])$, \ $\alpha\force F[\vec{a},0]$ \ and \ $\pi\in\|F[\vec{a},0]\|$. We must show that $I\star\ul{n}\ps\xi\ps\alpha\ps\pi\in\bbot$~;
this follows from lemma~\ref{lem_rec_eps}, with $k=0$.
\qed

\begin{lemma}\label{lem_rec_eps}
Let $n,k\in\NN$, \ $\xi\force\pt y(F[y]\to F[sy])$, \ $\alpha\force F[s^k0]$ \ and \ $\pi\in\|F[s^kn]\|$.\\
Then \ $\ul{n}\star\xi\ps\alpha\ps\pi\in\bbot$.
\end{lemma}

\proof
The proof is done for all integers $k$, by recurrence on $n$.
This is immediate if $n=0$.\\
In order to go from $n$ \ to $n+1$, we suppose now \ $\pi\in\|F[s^k(n+1)]\|$, \ i.e. \
$\pi\in\|F[s^{k+1}n]\|$.\\
We have \ $\ul{n+1}\star\xi\ps\alpha\ps\pi\succ\sig\ul{n}\star\xi\ps\alpha\ps\pi\succ
\sigma\star\ul{n}\ps\xi\ps\alpha\ps\pi\succ\ul{n}\star\xi\ps\xi\alpha\ps\pi$.\\
But, by hypothesis on $\xi$, we have \
$\xi\force F[s^k0]\to F[s^{k+1}0]$~; thus \ $\xi\alpha\force F[s^{k+1}0]$.\\
Hence the result, by the recurrence hypothesis.
\qed

\smallskip\noindent
{\bfseries Definition.} We denote by \ int$(n)$ \ the formula \
$\pt x(\pt y(sy\neps x\to y\neps x),n\neps x\to 0\neps x)$.

\smallskip\noindent
Theorem~\ref{storage} shows that the formulas \ int$(n)$ and $n\eps\wt{\NN}$ are interchangeable, i.e.
the formula \ $\pt n($int$(n)\dbfl n\eps\wt{\NN})$ is realized by a proof-like term~: this is the \emph{storage theorem for integers}.

\begin{lemma}\label{comp_succ}
$\lbd g\lbd x(g)(\sig)x\force\pt y(sy\neps\wt{\NN}\to y\neps\wt{\NN})$.
\end{lemma}

\proof
We show that \ $\lbd g\lbd x(g)(\sig)x\force sb\neps\wt{\NN}\to b\neps\wt{\NN}$ \ for every
individual $b$.\\
This is obvious if $b$ is not of the form $s^n0$, since then \ $\|b\neps\wt{\NN}\|=\vide$.
Thus, it remains to show~:\\
$\lbd g\lbd x(g)(\sig)x\force s^{n+1}0\neps\wt{\NN}\to s^n0\neps\wt{\NN}$. Thus, let
$\xi\force s^{n+1}0\neps\wt{\NN}$~; we must show~:\\
$\lbd g\lbd x(g)(\sig)x\star\xi\ps\ul{n}\ps\pi\in\bbot$, i.e.
$\xi\star\sig\ul{n}\ps\pi\in\bbot$, which is clear, since \ $\sig\ul{n}=\ul{n+1}$.
\qed

\begin{theorem}[Storage theorem]\label{storage}\ \\
\emph{\phantom ii)}~~$I\force\pt x^{\wt{\NN}}$ int$(x)$.\\
\emph{ii)}~~$T\force\pt x($int$(x),x\neps\wt{\NN}\to\bot)$ \ with \
$T=\lbd n\lbd f((n)\lbd g\lbd x(g)(\sig)x)f\ul{0}$.
\end{theorem}

\proof\hfill\\
\noindent\phantom ii)~~It is theorem~\ref{rec_eps}(i), if we take for $F[x,y]$ the formula \ $y\neps x$.

\smallskip\noindent
ii)~~Let \ $\nu\force$ int$(a)$, \ $\phi\force a\neps\wt{\NN}$ \ and \ $\pi\in\Pi$.
We must show \ $T\star\nu\ps\phi\ps\pi\in\bbot$, that is~:\\
$\nu\star\lbd g\lbd x(g)(\sig)x\ps\phi\ps\ul{0}\ps\pi\in\bbot$.\\
By hypothesis, we have \ $\nu\force\pt y(sy\neps\wt{\NN}\to y\neps\wt{\NN}),a\neps\wt{\NN}\to0\neps\wt{\NN}$.\\
But we have \ $\ul{0}\ps\pi\in\|0\neps\wt{\NN}\|$ \ by definition of \ $\wt{\NN}$ \ and, by lemma~\ref{comp_succ}~:\\
$\lbd g\lbd x(g)(\sig)x\force\pt y(sy\neps\wt{\NN}\to y\neps\wt{\NN})$.
Hence the result.
\qed

\smallskip\noindent
From theorem~\ref{rec_eps}(ii), it follows immediately that the \emph{recurrence scheme of ZF} is
realized in ${\mathcal N}$~; it is the scheme~:\\
$\pt\vec{x}(\pt y(F[\vec{x},y]\to F[\vec{x},sy]),F[\vec{x},0]\to(\pt n\in\wt{\NN})F[\vec{x},n])$ \
for every formula $F[\vec{x},y]$ of ZF (i.e. written with $\notin,\subseteq,0,s$).\\
Then, indeed, the formula $F$ is compatible with the extensional equivalence $\simeq$.

\smallskip\noindent
Since the function $s$ is compatible with $\simeq$, we deduce from lemma~\ref{comp_succ} that the
formula~:\\
$\pt y(y\in\wt{\NN}\to sy\in\wt{\NN})$ is realized in ${\mathcal N}$~; \ the formula \ $0\in\wt{\NN}$ \
is also obviously realized.\\
From the recurrence scheme just proved, we deduce that~:

\smallskip
\centerline{\emph{$\wt{\NN}$ is the set of integers of the model ${\mathcal N}$, considered as a
model of ZF\/.}}

\begin{theorem}\label{recursive}\ \\
\emph{\phantom ii)}~~Let $f:\NN^k\to\NN$ be a recursive function. Then, the formula~:\\
$\pt x_1^{\wt{\NN}}\ldots\pt x_k^{\wt{\NN}}(f(x_1,\ldots,x_k)\eps\wt{\NN})$ \
is realized in ${\mathcal N}$.\\
\emph{ii)}~~Let $g:\NN^k\to2$ be a recursive function. Then, the formula~:\\
$\pt x_1^{\wt{\NN}}\ldots\pt x_k^{\wt{\NN}}(g(x_1,\ldots,x_k)=1\lor
g(x_1,\ldots,x_k)=0)$ \ is realized in ${\mathcal N}$.
\end{theorem}\noindent
i)~~This can be written \ $\pt x_1\inde\ldots\pt x_k\inde$ ent$(f(x_1,\ldots,x_k))$. \
The proof is done in \cite{krivine6,krivine3}.\\
ii)~~We have \
${\mathcal N}\force(\pt x_1\eps\gn)\ldots(\pt x_k\eps\gn)\,
g(x_1,\ldots,x_k)\eps\gd$.\\
Now, since $g$ is recursive, we have, by (i)~:\\
${\mathcal N}\force(\pt x_1\eps\wt{\NN})\ldots(\pt x_k\eps\wt{\NN})\,
g(x_1,\ldots,x_k)\eps\wt{\NN}$.\\
Hence the result, by lemma~\ref{bool}.
\qed

\begin{lemma}\label{bool}
$\lbd x\lbd y\lbd f(f)xy\force\pt x^{\gd}(x\ne1,x\ne0\to x\neps\wt{\NN})$.
\end{lemma}

\proof
We have to show~:\\
$\lbd x\lbd y\lbd f(f)xy\force\top,\bot\to 0\neps\wt{\NN}$ \ and \
$\lbd x\lbd y\lbd f(f)fxy\force\bot,\top\to 1\neps\wt{\NN}$.\\
Thus let \ $\xi\force\top$ (i.e. $\xi\in\Lbd$ arbitrary) and \ $\eta\force\bot$.
We have to show~:\\
$\lbd x\lbd y\lbd f(f)xy\star\xi\ps\eta\ps\ul{0}\ps\pi\in\bbot$ \ and \
$\lbd x\lbd y\lbd f(f)xy\star\eta\ps\xi\ps\ul{1}\ps\pi\in\bbot$\\
which is trivial.
\qed

\smallskip\noindent
{\small{\bfseries Remarks.} i)~~In the present paper, theorem~\ref{recursive} \ is used only
in trivial particular cases.\\
ii)~~Let us recall the difference between $\gn$ and $\wt{\NN}$ (the set of integers in the
model ${\mathcal N}$)~; we have~:\\
$\xi\force\pt x^{\gn}\,F[x]$ \ iff \
$(\pt n\in\NN)(\pt\pi\in\|F[s^n0]\|)\,\xi\star\pi\in\bbot$.\\
$\xi\force\pt x^{\wt{\NN}}\,F[x]$ \ iff \
$(\pt n\in\NN)(\pt\pi\in\|F[s^n0]\|)\,\xi\star\ul{n}\ps\pi\in\bbot$.\\
Notice that we have \ $K\force\pt x(x\neps\gn\to x\neps\wt{\NN})$,
in other words \ $K\force\wt{\NN}\subset\gn$. This means that, in ${\mathcal N}$,
the set $\wt{\NN}$ of integers is strongly included in $\gn$. In the particular
realizability model considered below (and, in fact, in every non trivial realizability
model), the formula \ $\gn\not\subseteq\wt{\NN}$ \ is realized.}

\subsection*{Non extensional and dependent choice}\ \\
For each formula $F(x,y_1,\ldots,y_m)$ of \ZFe, we add a function symbol $f_F$ of arity $m+1$, with the axiom~: \ \
$\pt\vec{y}(\pt k^{\wt{\NN}}F[f_F(k,\vec{y}),\vec{y}]\to\pt x\,F[x,\vec{y}])$\\
or else~: \ \ $\pt\vec{y}(\pt k\inde F[f_F(k,\vec{y}),\vec{y}]\to\pt x\,F[x,\vec{y}])$.

\smallskip\noindent
It is the {\em axiom scheme of non extensional choice}, in abbreviated form NEAC.

\smallskip\noindent
{\small{\bfseries Remarks.}
i)~~The axiom scheme NEAC does not imply the axiom of choice in ZF\/, because we do not suppose that the
symbol $f_F$ is compatible with the extensional equivalence $\simeq$. It is the reason why we speak about
\emph{non extensional} axiom of choice. On the other hand, as we show below, it implies DC (the axiom of dependent choice).\\
ii)~~It seems that we could take for $f_F$ a $m$-ary function symbol and use the following simpler
(and logically equivalent) axiom scheme NEAC'~: \
$\pt\vec{y}(F[f_F(\vec{y}),\vec{y}]\to\pt x\,F[x,\vec{y}])$.\\
But this axiom scheme cannot be realized, even though the axiom scheme NEAC is realized by a
very simple proof-like term (theorem~\ref{ACNE}), \emph{provided the instruction $\vsig$ is present}.\\
More precisely, we can define a function $f_F$ in ${\mathcal M}$, such that NEAC is realized
in ${\mathcal N}$, but this is impossible for NEAC'.}

\begin{theorem}[NEAC]\label{ACNE}\ \\
For each closed formula \ $\pt x\pt\vec{y}\,F$, we can define a $(m+1)$-ary function
symbol $f_F$ such that~:\\
$\lbd x(\vsig)xx\force\pt\vec{y}(\pt k\inde F[f_F(k,\vec{y})/x,\vec{y}]\to
\pt x\,F[x,\vec{y}])$.
\end{theorem}

\proof
For each $k\in\ennl$ we put \ $P_k=\{\pi\in\Pi;$ $\xi\star\ul{k}\ps\pi\notin\bbot,\,k=\nn_\xi\}$.

\smallskip\noindent
For each individual $x$, we have~: \ $\dsp\|\pt x\,F[x,\vec{y}]\|=\bigcup_a\|F[a,\vec{y}]\|$.\\
Therefore, there exists a function $f_F$ such that, given $k\in\ennl$ and $\vec{y}$ such that\\
$P_k\cap\|\pt x\,F[x,\vec{y}]\|\ne\vide$,  \ we have \ $P_k\cap\|F[f_F(k,\vec{y}),\vec{y}]\|\ne\vide$.

\smallskip\noindent
Now, we want to show $\lbd x(\vsig)xx\force\pt k\inde F[f_F(k,\vec{y}),\vec{y}]\to F[x,\vec{y}]$,
for every individuals $x,\vec{y}$.\\
Thus, let \ $\xi\force\pt k\inde F[f_F(k,\vec{y}),\vec{y}]$ and
$\pi\in\|F[a,\vec{y}]\|$~; we must show \ $\lbd x(\vsig)xx\star\xi\ps\pi\in\bbot$.\\
If this is false, we have \ $\vsig\star\xi\ps\xi\ps\pi\notin\bbot$ \ and therefore
$\xi\star\ul{j}\ps\pi\notin\bbot$ with $j=\nn_\xi$.\\
It follows that \ $\pi\in P_j\cap\|F[a,\vec{y}]\|$~; thus, there exists \
$\pi'\in P_j\cap\|F[f_F(j,\vec{y}),\vec{y}]\|$.\\
Now, we have \ $\ul{j}\ps\pi'\in\|\pt k\inde F[f_F(k,\vec{y}),\vec{y}]\|$, and therefore,
by hypothesis on $\xi$, we have~:\\
$\xi\star\ul{j}\ps\pi'\in\bbot$.  This is in contradiction with $\pi'\in P_j$.
\qed

\subsubsection*{NEAC implies DC}
Let us call DCS (dependent choice scheme) the following axiom scheme~:\\
$\pt\vec{z}(\pt x\ex y\,F[x,y,\vec{z}]\!\to\!\pt n\inde\ex!y\,S_F[n,y,\vec{z}]\land
\pt n\inde\ex y\ex y'\{S_F[n,y,\vec{z}],S_F[sn,y',\vec{z}],F[y,y',\vec{z}]\})$.\\
where $F$ is a formula of \ZFe\ with free variables $x,y,\vec{z}$~; the formula $S_F$ is written
below.\\
In the following, we omit the variables $\vec{z}$ (the parameters), for sake of simplicity.\\
The usual axiom of dependent choice DC is obtained by taking for $F[x,y,z_0,z_1]$ the formula \
$y\eps z_0\land(x\eps z_0\to{<}x,y{>}\eps z_1)$.\\
We now show how to define the formula $S_F$, so that \ \ZFe, NEAC $\vdash$ DCS~; \
we shall conclude that DC is realized.

\smallskip\noindent
So, let us assume $\pt x\ex y\,F[x,y]$. By NEAC, there is a function symbol $f$ such that~:\\
$\pt x\ex k\inde F[x,f(k,x)]$. We define the formula \ $R_F[x,y]$ as follows~:\\
$R_F[x,y]\equiv\ex k\inde\{F[x,f(k,x)],\pt i\inde(i<k\to\neg F[x,f(i,x)]),y=f(k,x)\}$.\\
This means~: ``$y=f(k,x)$ for the first integer $k$ such that $F[x,f(k,x)]$~''.\\
Therefore, $R_F$ is functional, i.e. we have \ $\pt x\ex!y\,R_F(x,y)$.\\
$S_F$ is defined so as to represent a sequence obtained by iteration of the function given by $R_F$,
beginning (arbitrarily) at $0$~:\\
$S_F(n,x)\equiv\pt z[\pt m\pt y\pt y'({<}m,y{>}\eps z,R_F(y,y')\to{<}sm,y'{>}\eps z),{<}0,0{>}\eps z
\to{<}n,x{>}\eps z]$.\\
It should be clear that, with this definition of $S_F$, we obtain~:\\
$\pt n\inde\ex!y\,S_F[n,y]$ \ and \ $\pt n\inde\ex y\ex y'\{S_F[n,y],S_F[sn,y'],F[y,y']\}$.\\
Thus, DCS is provable from \ZFe \ and NEAC.

\smallskip\noindent
{\small{\bfseries Remark.} We have used the binary function symbol ${<}x,y{>}$ which is defined,
in the ground model~${\mathcal M}$, in the usual way~: ${<}a,b{>}=\{\{a\},\{a,b\}\}$. Then, the formulas~:\\
$\pt x\pt x'\pt y\pt y'({<}x,y{>}={<}x',y'{>}\mapto x=x')$, \
$\pt x\pt x'\pt y\pt y'({<}x,y{>}={<}x',y'{>}\mapto y=y')$,\\
are trivially realized by $I$.}

\subsection*{Properties of the Boolean algebra $\gd$}\ \\
Let $(x\ppt y)$ be the binary recursive function  defined as follows in ${\mathcal M}$~:\\
$(m\ppt n)=1$ if $m,n\in\NN,\,m<n$~; \ else $(m\ppt n)=0$.

\begin{theorem}\label{ord_ent}
For every choice of \ $\bbot$, the relation $(x\ppt y)=1$ is, in ${\mathcal N}$, a strict well
founded partial order, which is the usual order on integers (i.e. on $\wt{\NN}$).
\end{theorem}

\proof
Indeed, the formulas~:\\
$\pt x((x\ppt x)\ne1)$ and $\pt x\pt y\pt z((x\ppt y)=1\mapto((y\ppt z)=1\mapto (x\ppt z)=1))$\\
are trivially realized.\\
Moreover, since the relation $(x\ppt y)=1$ is well founded, we have (theorem~\ref{bien_fonde})~:\\
$\Y\force\pt x(\pt y((y\ppt x)=1\mapto F[y])\to F[x])\to\pt x\,F[x]$\\
for every formula \ $F[x]$ with parameters and one free variable.

\smallskip\noindent
By theorem~\ref{recursive}(ii), the binary recursive function $(x\ppt y)$ sends $\wt{\NN}^2$ into $\{0,1\}$, in the model~${\mathcal N}$.
Therefore, it suffices to check that the following formulas are realized in ${\mathcal N}$~:\\
$\pt x^{\wt{\NN}}\pt y^{\wt{\NN}}(y\le x\to(x\ppt y)\ne1)$~; \
$\pt x^{\wt{\NN}}\pt y^{\wt{\NN}}(x<y\to(x\ppt y)=1)$.\\
Now the following formulas are trivially realized~:\\
$\pt x^{\gn}\pt y^{\gn}\pt z^{\gn}(x=y+z\to(x\ppt y)\ne1)$~; \
$\pt x^{\gn}\pt y^{\gn}\pt z^{\gn}(y=x+z+1\to(x\ppt y)=1)$.
\qed

\smallskip\noindent
In the ground model ${\mathcal M}$, we put, for each integer $n$~:\\
\centerline{$\mb{n}=\{0,1,\ldots,n-1\}=\{0,s0,\ldots,s^{n-1}0\}$.}

\noindent
The functions $n\longmapsto\mb{n}$ and $n\longmapsto\gimel\mb{n}$ are defined
in the realizability model ${\mathcal N}$, with domain~$\gimel\NN$.

\begin{theorem}\label{restriction_finie}\ \\
The following formulas are realized in ${\mathcal N}$~:\\
\emph{\phantom{ii}i)}~~$\pt x^{\gn}\pt m^{\gn}((x\ppt m)=1\dbfl x\eps\gimel\mb{m})$~;\\
\emph{\phantom iii)}~~$\pt m^{\gn}\pt n^{\gn}((m\ppt n)=1\to\gimel\mb{m}\subset\gimel\mb{n})$~;\\
\emph{iii)}~~$\pt x^{\gn}\pt m^{\gn}((x\ppt m)=1\dbfl\ex y^{\gn}(m=x+y+1))$.
\end{theorem}

\proof
Remember that \ $x\subset y$ is the formula \ $\pt z(z\neps y\to z\neps x)$.

\smallskip\noindent
\phantom{ii}i)~~We have trivially \ $\|(a\ppt m)\ne1\|=\|a\neps\gimel\mb{m}\|$ for every $a,m\in\NN$.\\
\phantom iii)~~By transitivity of the relation $(m\ppt n)=1$ (theorem~\ref{ord_ent}).\\
iii)~~We observe that \ $\|(a\ppt m)\ne1\|=\|(\pt y\eps\gn)(m\ne a+y+1)\|$ \ for every \ $a,m\in\NN$.
\qed

\smallskip\noindent
For each \ $n\eps\gn$ (and, in particular, for each \ $n\eps\wt{\NN}$, i.e.~for each integer of
${\mathcal N}$), the set defined, in ${\mathcal N}$, by $(x\ppt n)=1$ (the strict initial segment defined by $n$)
is therefore extensionally equivalent to $\gimel\mb{n}$.

\begin{theorem}\label{B_mn}
In ${\mathcal N}$, the application $(x,y)\longmapsto my+x$ is a bijection from $\gimel\mb{m}\fois\gimel\mb{n}$
onto $\gimel(\mb{mn})$. Indeed, the following formulas are realized in ${\mathcal N}$ by $I$~:\\
\emph{\phantom{ii}i)}~~$\pt m^{\gn}\pt n^{\gn}\pt x^{\gimel\mb{m}}\pt y^{\gimel\mb{n}}((my+x)\eps\gimel\mb{mn})$~;\\
\emph{\phantom iii)}~~$\pt m^{\gn}\pt n^{\gn}\pt x^{\gimel\mb{m}}\pt x'^{\gimel\mb{m}}
\pt y^{\gimel\mb{n}}\pt y'^{\gimel\mb{n}}(my+x=my'+x'\mapto x=x')$~;\\
\hspace*{1.8em}$\pt m^{\gn}\pt n^{\gn}\pt x^{\gimel\mb{m}}\pt x'^{\gimel\mb{m}}
\pt y^{\gimel\mb{n}}\pt y'^{\gimel\mb{n}}(my+x=my'+x'\mapto y=y')$~;\\
\emph{iii)}~~$\pt m^{\gn}\pt n^{\gn}\pt z^{\gimel\mb{mn}}\ex x^{\gimel\mb{m}}\ex y^{\gimel\mb{n}}(z=my+x)$.
\end{theorem}

\proof\hfill\\
\noindent\phantom{ii}i)~and~\phantom iii)~~We simply have to replace $\pt m^{\gn}$ \ and \ $\pt x^{\gimel\mb{m}}$ \ with
their definitions, which are~: \
$\pt m^{\gn}F\equiv\pt m(1_{\NN}(m)=1\mapto F)$\ ; \
$\pt x^{\gimel\mb{m}}F\equiv\pt x((x\ppt m)=1\mapto F)$.\\
We see immediately that these two formulas are realized by $I$.

\smallskip\noindent
iii)~~We show that~:\\
$I\force\pt m^{\gn}\pt n^{\gn}\pt z^{\gn}(\pt x^{\gn}\pt y^{\gn}((x\ppt m)=1\mapto
((y\ppt n)=1\mapto z\ne my+x))\to(z\ppt mn)\ne 1)$.\\
Thus, we consider~:\\
$m,n,z_0\in\NN$~; \ $\xi\in\Lbd,\,\xi\force\pt x^{\gn}\pt y^{\gn}((x\ppt m)=1\mapto((y\ppt n)=1\mapto z\ne my+x))$\\
and \ $\pi\in\|(z_0\ppt mn)\ne 1\|$. We must show \ $I\star\xi\ps\pi\in\bbot$, that is $\xi\star\pi\in\bbot$.\\
We have \ $\|(z_0\ppt mn)\ne 1\|\ne\vide$, \ therefore \ $z_0<mn$.\\
Thus, there exist $x_0,y_0\in\NN,x_0<m,y_0<n$ \ such that \ $z_0=mx_0+y_0$.\\
Now, by hypothesis on $\xi$, we have~:\\
$\xi\force(x_0\ppt m)=1\mapto((y_0\ppt n)=1\mapto z_0\ne my_0+x_0)$, in other words \ $\xi\force\bot$.
\qed

\subsubsection*{Injection of $\gimel\mb{n}$ into ${\mathcal P}(\wt{\NN})$}
Remember that we have fixed a recursive bijection~: $\xi\longmapsto\nn_\xi$ \ from $\Lbd$
onto \ $\NN$.
The inverse bijection will be denoted \ $n\longmapsto\xi_n$.\\
This bijection is used in the execution rule of the instruction $\vsig$, which is as follows~:\\
\centerline{$\vsig\star\xi\ps\eta\ps\pi\succ\xi\star\ul{\nn}\,_\eta\ps\pi$.}

\smallskip\noindent
We define, in ${\mathcal M}$, a function \ $\Delta:\NN\to2$ \ by putting \
$\Delta(n)=0$ $\Dbfl$ $\xi_n\force\bot$.\\
In this way, we have defined a function symbol $\Delta$, in the language of \ZFe. In the
realizability model~${\mathcal N}$, the symbol $\Delta$ represents a function from $\gn$
into $\gd$. In particular, the function $\Delta$ sends the set $\wt{\NN}$ of
integers of the model ${\mathcal N}$ into the Boolean algebra $\gd$.

\begin{theorem}\label{vsig_yxx}
Let us put \ $\theta=\lbd x\lbd y(\vsig)yxx$~; then, we have~:\\
\centerline{
$\theta\force\pt x^{\gd}(x\ne0\to\ex n\inde\{\Delta(n)\ne0,\Delta(n)\le x\})$}
where $\le$ is the order relation of the Boolean algebra \ $\gd$~: \ $y\le x$ \ is the
formula \ $x=(y\,\ou\,x)$.
\end{theorem}

\proof
We must show \
$\theta\force\pt x^{\gd}(x\ne0,\pt n\inde(\Delta(n)\ne0\to
x\ne\Delta(n)\ou\,x)\to\bot)$.\\
Thus, let \ $a\in\{0,1\}$, \ $\xi\force a\ne0$, \
$\eta\force\pt n\inde(\Delta(n)\ne0\to a\ne\Delta(n)\ou\,a)$ \ and \ $\pi\in\Pi$.\\
We must show \ $\theta\star\xi\ps\eta\ps\pi\in\bbot$ \ that is \ $\vsig\star\eta\ps\xi\ps\xi\ps\pi\in\bbot$, or else \
$\eta\star\ul{\nn}_\xi\ps\xi\ps\pi\in\bbot$.\\
By hypothesis on $\eta$, it suffices to show \
$\ul{\nn}_\xi\ps\xi\ps\pi\in\|\pt n\inde(\Delta(n)\ne0\to a\ne\Delta(n)\ou\,a)\|$, \
that is, by definition of the quantifier \ $\pt n\inde$~:\hspace{1em}
$\xi\ps\pi\in\|\Delta(\nn_\xi)\ne0\to a\ne\Delta(\nn_\xi)\ou\,a\|$.\\
This amounts to show \ $\xi\force\Delta(\nn_\xi)\ne0$ \ and \ $a=\Delta(\nn_\xi)\ou\,a$.

\smallskip\noindent
$\bullet$~~Proof of \ $\xi\force\Delta(\nn_\xi)\ne0$~: \
if $\Delta(\nn_\xi)=1$, this is trivial, because \ $\|\Delta(\nn_\xi)\ne0\|=\vide$~;\\
if \ $\Delta(\nn_\xi)=0$, \ then $\xi\force\bot$, \ by definition of $\Delta$.\\
$\bullet$~~Proof of \ $a=\Delta(\nn_\xi)\ou a$~: \ this is obvious if $a=1$~; if $a=0$, then \ $\xi\force\bot$, \ by hypothesis on $\xi$.
Therefore $\Delta(\nn_\xi)=0$ by definition of $\Delta$, hence the result.
\qed

\smallskip\noindent
By theorem~\ref{vsig_yxx}, the set $\{\Delta(n);\;n\eps\wt{\NN},\Delta(n)\ne0\}$ is,
in the realizability model~${\mathcal N}$, a~countable dense subset of the Boolean algebra
$\gd$~: this means that each element $\ne0$ of this Boolean algebra
has a lower bound of the form $\Delta(n)$, with $n\eps\wt{\NN}$ and $\Delta(n)\ne0$.\\
It follows that the application of \ $\gd$ \ into \ ${\mathcal P}(\wt{\NN})$ given by~:

\centerline{$x\longmapsto\{n\eps\wt{\NN};\;\Delta(n)\le x,\Delta(n)\ne0\}$}

\noindent
is one to one~: indeed, if \ $a,b\eps\gd$ with $a\ne b$, then $a+b\ne0$~; thus, there
exists an integer $n\eps\wt{\NN}$ such that $\Delta(n)\ne0$ and $\Delta(n)\le a+b$.
Therefore, we have \ \ $\Delta(n)\le a$ \ \ iff \ \ $(b\et\Delta(n))=0$.\\
But, since \ $\Delta(n)\ne0$, we get~: \ \
$\Delta(n)\le a$ \ \ iff \ \ $\Delta(n)\not\le b$.\\
We have shown~:

\begin{theorem}\label{gimel2_R}\ \\
The formula~: ``there exists an injection of $\gd$ into ${\mathcal P}(\wt{\NN})$'' is realized in the model ${\mathcal N}$.\qed
\end{theorem}

\begin{corollary}\label{gimel_n_R}
The formula~: ``for every integer $n$ there exists an injection of $\gimel\mb{n}$
into ${\mathcal P}(\wt{\NN})$'' is realized in the model ${\mathcal N}$.
\end{corollary}

\proof
Using theorem~\ref{B_mn} we see, by recurrence on $m$, that the model~${\mathcal N}$ realizes the formula~:\\
``~$\pt m^{\wt{\NN}}((\gd)^m$ is equipotent to $\gimel(\mb{2^m}))$~''~; \ \ and therefore
also the formula~:\\
``~$\pt m^{\wt{\NN}}($there exists an injection of \ $\gimel(\mb{2^m})$ \ into \
${\mathcal P}(\wt{\NN}))$~''.\\
Finally, by theorem~\ref{restriction_finie}(ii), we see that the following formula is realized~:\\
``~$\pt n^{\wt{\NN}}($there exists an injection of \ $\gimel\mb{n}$ \ into \
${\mathcal P}(\wt{\NN}))$~''.
\qed

\section{Realizability models in which $\mathbb{R}$ is not well ordered}

\subsection*{$\gd$ atomless}

\begin{theorem}\label{B_non_triviale}
We suppose there exist two proof-like terms \ $\omega_0,\omega_1$ \ such that, for every \
$\pi\in\Pi$, we have \ $\omega_0\kk_\pi\force\bot$ or \ $\omega_1\kk_\pi\force\bot$. Then, the Boolean algebra $\gd$ is non trivial. Indeed~:
$\theta\force\pt x(x\ne1,x\ne0\to x\neps\gd)\to\bot$ \ with \
$\theta=\lbd f(\Ccc)\lbd k((f)(\omega_1)k)(\omega_0)k$.
\end{theorem}

\proof
Let \ $\xi\force\pt x(x\ne1,x\ne0\to x\neps\gd)$ \ and \ $\pi\in\Pi$. We must show~:\\
$\theta\star\xi\ps\pi\in\bbot$, \ that is \ $\xi\star\omega_1\kk_\pi\ps\omega_0\kk_\pi\ps\pi\in\bbot$.\\
But, by hypothesis on $\xi$, we have \ $\xi\force\top,\bot\to\bot$ \ and \ $\xi\force\bot,\top\to\bot$. Hence the result, by hypothesis on \ $\omega_1,\omega_0$.
\qed

\smallskip\noindent
{\small{\bfseries Remark.} When the Boolean algebra $\gd$ is non trivial, there are
necessarily non standard integers in the realizability model ${\mathcal N}$, i.e. integers
which are not in ${\mathcal M}$. Indeed, let $a\eps\gd,a\ne0,1$~; by theorem~\ref{vsig_yxx},
there is an integer $n$ such that $\Delta(n)\ne0,\Delta(n)\le a$~; thus $\Delta(n)\ne1$.
The integer $n$ cannot be standard, since $\Delta(m)=0$ or $1$ if $m$ is in ${\mathcal M}$.}

\begin{theorem}\label{sans_atome}
We suppose that there exists three proof-like terms \ $\alpha_0,\alpha_1,\alpha_2$ \
such that, for every \ $\xi\in\Lbd$ and $\pi\in\Pi$, we have \ $\kk_\pi\xi\alpha_0\force\bot$ or \
$\kk_\pi\xi\alpha_1\force\bot$ \ or \ $\kk_\pi\xi\alpha_2\force\bot$.\\
Then, the Boolean algebra $\gd$ is atomless. Indeed~:\\
$\theta\force\pt x[\pt y(x\et y\ne0,x\et y\ne x\to y\neps\gd),x\ne0\to x\neps\gd]$\\
with \ $\theta=\lbd x\lbd y(\Ccc)\lbd k((x)(k)y\alpha_0)((x)(k)y\alpha_1)(k)y\alpha_2$.
\end{theorem}

\proof
By a simple computation, we see that we must show~:\\
\phantom ii) $\theta\force(\bot,\bot\to\bot),\bot\to\bot$.\\
ii) $\theta\force|\top,\bot\to\bot|\cap|\bot,\top\to\bot|,\top\to\bot$.\\
Proof of (i)~: let \ $\eta\in|\bot,\bot\to\bot|$ \ and \ $\xi\in|\bot|$.
We must show \ $\theta\star\eta\ps\xi\ps\pi\in\bbot$, that is~:\\
$\eta\star\kk_\pi\xi\alpha_0\ps((\eta)(\kk_\pi)\xi\alpha_1)(\kk_\pi)\xi\alpha_2\ps\pi\in\bbot$.\\
But, from \ $\xi\force\bot$, we deduce \ $\kk_\pi\xi\zeta\force\bot$ for every $\zeta\in\Lambda_c$.\\
Since \ $\eta\force\bot,\bot\to\bot$, \ we have \ $((\eta)(\kk_\pi)\xi\alpha_1)(\kk_\pi)\xi\alpha_2\force\bot$ \ and therefore~:\\
$\eta\star\kk_\pi\xi\alpha_0\ps((\eta)(\kk_\pi)\xi\alpha_1)(\kk_\pi)\xi\alpha_2\ps\pi\in\bbot$.

\smallskip\noindent
Proof of (ii)~: let \ $\eta\in|\top,\bot\to\bot|\cap|\bot,\top\to\bot|$ \ and \
$\xi\in\Lambda_c$. Again, we must show that~:\\
$\eta\star\kk_\pi\xi\alpha_0\ps((\eta)(\kk_\pi)\xi\alpha_1)(\kk_\pi)\xi\alpha_2\ps\pi\in\bbot$. If this is false, then~:\\
$\kk_\pi\xi\alpha_0\nforce\bot$ (because \ $\eta\force\bot,\top\to\bot$) and\\
$((\eta)(\kk_\pi)\xi\alpha_1)(\kk_\pi)\xi\alpha_2\nforce\bot$ (because $\eta\force\top,\bot\to\bot$).\\
But, since \ $\eta\force\bot,\top\to\bot$ (resp. $\top,\bot\to\bot$), we have \
$\kk_\pi\xi\alpha_1\nforce\bot$ (resp. $\kk_\pi\xi\alpha_2\nforce\bot$).\\
This contradicts the hypothesis of the theorem.
\qed

\subsection*{$\mathbb{R}$ not well orderable}

\begin{theorem}\label{pas_surj}\ \\
We suppose that there exists a proof-like term \ $\omega$ such that, for every \ $\xi,\xi'\in\Lbd,\,\xi\ne\xi'$ and
$\pi\in\Pi$, we have \ $\omega\kk_\pi\xi\force\bot$ \ or \ $\omega\kk_\pi\xi'\force\bot$.\\
Then we have, for every formula $F$ with three free variables~:\\
$\theta\force\pt m^{\gn}\pt n^{\gn}\pt z[(m\ppt n)=1\mapto\\
\hspace*{\fill}(\pt x\pt y\pt y'(F(x,y,z),F(x,y',z),y\ne y'\to\bot),
\pt y^{\gimel\mb{n}}\neg\pt x^{\gimel\mb{m}}\neg\,F(x,y,z)\to\bot)]$\\
with \ $\theta=\lbd x\lbd x'(\Ccc)\lbd k(x')\lbd z(xzz)(\omega)kz$.
\end{theorem}

\noindent{\small{\bfseries Remark.}
This shows that, if $(m\ppt n)=1$, then ($\gimel\mb{m}\subset\gimel\mb{n}$ and) there is no surjection of
$\gimel\mb{m}$ onto~$\gimel\mb{n}$~: indeed, it suffices to take, for $F(x,y,z)$, the formula ${<}x,y{>}\eps z$.}

\smallskip\proof
Assume this is false~; then, there exist $m,n\in\ennl$ with $m<n$, an individual $c$, two terms $\xi,\xi'\in\Lbd$ and a stack $\pi\in\Pi$ such that~:\\
$\theta\star\xi\ps\xi'\ps\pi\notin\bbot$~;\\
$\xi\force\pt x\pt y\pt y'[F(x,y,c),F(x,y',c),y\ne y'\to\bot]$~;\\
$\xi'\force\pt y^{\gimel\mb{n}}\neg\pt x^{\gimel\mb{m}}\neg\,F(x,y,c)$.

\smallskip\noindent
Therefore, we have \ $\xi'\star\eta\ps\pi\notin\bbot$ \ with \ $\eta=\lbd z(\xi zz)(\omega)\kk_\pi z$.
By hypothesis on \ $\xi'$ \ we have, for every integer \ $i<n$~: \
$\eta\nforce\pt x^{\gimel\mb{m}}\neg F(x,i,c)$. Thus, there exists an integer \ $m_i<m$ \ such that
$\eta\nforce\neg F(m_i,i,c)$. It follows that there exist \ $\xi_i\in\Lbd$ \ and \ $\pi_i\in\Pi$ \
such that $\xi_i\force F(m_i,i,c)$ \ and \ $\eta\star\xi_i\ps\pi_i\notin\bbot$. By definition of $\eta$,
we get $\xi\star\xi_i\ps\xi_i\ps\omega\kk_\pi\xi_i\ps\pi_i\notin\bbot$. By hypothesis
on \ $\xi$, it follows that \ $\omega\kk_\pi\xi_i\nforce i\ne i$~; in other words, we have \ $\omega\kk_\pi\xi_i\nforce\bot$ \ for every integer \ $i<n$.\\
By the hypothesis of the theorem, it follows that we have \ $\xi_i=\xi_j$ \ for every \ $i,j<n$.\\
But, since \ $m_i<m<n$ \ and \ $i<n$, \ there exist \ $i,j<n$, $i\ne j$ \ such that \ $m_i=m_j=k$.\\
Then, \ $\xi_i=\xi_j\force F(k,i,c),F(k,j,c)$ \ and
$\omega\kk_\pi\xi_i\force i\ne j$ since $\|i\ne j\|=\vide$.\\
Therefore, by hypothesis on $\xi$, we have \
$\xi\star\xi_i\ps\xi_i\ps\omega\kk_\pi\xi_i\ps\pi_i\in\bbot$, which is a contradiction.
\qed

\smallskip\noindent
Now, we see that, with the hypothesis of theorem~\ref{pas_surj}, there is no surjection from
$\gd$ onto\\
$\gd\fois\gd$. Indeed, by theorem~\ref{B_mn}, there exists a bijection from $\gd\fois\gd$
onto $\gimel\mb{4}$ and, by theorem~\ref{pas_surj}, there is no surjection from $\gd$ onto
$\gimel\mb{4}$. But, by theorem~\ref{sans_atome}, $\gd$ is infinite~; it follows that
{\em $\gd$ cannot be well ordered}.\\
Now\/, by theorem~\ref{gimel2_R}, $\gd$ is equipotent with a subset of \ ${\mathcal P}(\wt{\NN})$.
Therefore, the hypothesis of theorems~\ref{sans_atome} and~\ref {pas_surj} are sufficient in
order that the following formula be realized in the model~${\mathcal N}$~:\\
\centerline{\em There is no well ordering on the set of reals.}

\smallskip\noindent
In fact, the hypothesis of theorem~\ref{pas_surj} is sufficient~: this follows from
theorem~\ref{B_non_den}.

\begin{theorem}\label{B_non_den}
Same hypothesis as theorem~\ref{pas_surj}~: there exists a proof-like term $\omega$ such that, for every \
$\pi\in\Pi$ \ and \ $\xi,\xi'\in\Lbd,\,\xi\ne\xi'$, we have \
$\omega\kk_\pi\xi\force\bot$ \ or \ $\omega\kk_\pi\xi'\force\bot$.\\
Then we have, for every formula $F$ with three free variables~:\\
$\theta\force\pt z\{\pt x[\pt n\inde\,F(n,x,z)\to x\neps\gd],
\pt n\pt x\pt y[\neg F(n,x,z)\neg F(n,y,z),x\ne y\to\bot]\to\bot\}$\\
with \ $\theta=\lbd x\lbd x'(\Ccc)\lbd k(x)\lbd n(\Ccc)\lbd h(x'hh)(\omega k)\lbd f(f)hn$.
\end{theorem}\noindent
{\small{\bfseries Remark.}
This formula means that, in the realizability  model ${\mathcal N}$, there is no surjection
from the set of integers \ $\wt{\NN}$ \ onto \ $\gd$~: it suffices to take for
$F(x,y,z)$ the formula ${<}x,y{>}\neps z$ (the  graph of an hypothetical surjection
being \ ${<}x,y{>}\eps z$).}

\smallskip\proof
Reasoning by contradiction, we suppose that there is an individual $c$, a stack
$\pi\in\Pi$, and two terms $\xi,\xi'$ such that~:\\
$\xi\force\pt x[\pt n\inde\,F(n,x,c)\to x\neps\gd]$~; \ 
$\xi'\force\pt n\pt x\pt y[\neg F(n,x,c)\neg F(n,y,c),x\ne y\to\bot]$ \ and\\ $\theta\star\xi\ps\xi'\ps\pi\notin\bbot$.\\
Therefore, we have \ $\xi\star\eta\ps\pi\notin\bbot$, with \
$\eta=\lbd n(\Ccc)\lbd h(\xi'hh)(\omega\kk_\pi)\lbd f(f)hn$.\\
By hypothesis on $\xi$, we have \ $\eta\nforce\pt n\inde\,F(n,0,c)$ \ and \
$\eta\nforce\pt n\inde\,F(n,1,c)$. Thus, we see that there exist $n_0,n_1\in\ennl$,
$\pi_0\in\|F(n_0,0,c)\|$ \ and \ $\pi_1\in\|F(n_1,1,c)\|$ \ such that \ $\eta\star\ul{n}_0\ps\pi_0\notin\bbot$ \
and \ $\eta\star\ul{n}_1\ps\pi_1\notin\bbot$. By performing these two processes, we obtain~:\\ $\xi'\star\kk_{\pi_0}\ps\kk_{\pi_0}\ps\zeta_0\ps\pi_0\notin\bbot$ et
$\xi'\star\kk_{\pi_1}\ps\kk_{\pi_1}\ps\zeta_1\ps\pi_1\notin\bbot$,\\
with \ $\zeta_0=(\omega\kk_\pi)\lbd f(f)\kk_{\pi_0}\ul{n}_0$ \ and \
$\zeta_1=(\omega\kk_\pi)\lbd f(f)\kk_{\pi_1}\ul{n}_1$.\\
By hypothesis on $\xi'$, we have $\xi'\force\neg F(n_0,0,c),\neg F(n_0,0,c),0\ne 0\to\bot$.\\
Since \ $\kk_{\pi_0}\force\neg F(n_0,0,c)$, \ we see that \ $\zeta_0\nforce\bot$ \ and,
in the same way, \ $\zeta_1\nforce\bot$.\\
Thus, by the hypothesis of the theorem, we have~:\\
$\lbd f(f)\kk_{\pi_0}\ul{n}_0=\lbd f(f)\kk_{\pi_1}\ul{n}_1$, and therefore \ $n_0=n_1$ \ and \ $\pi_0=\pi_1$.\\
But, we have \ $\xi'\force\neg F(n_0,0,c),\neg F(n_0,1,c),0\ne 1\to\bot$. Moreover, we have~:\\
$\pi_0\in\|F(n_0,0,c)\|$ \ and \ $\pi_1\in\|F(n_1,1,c)\|$, thus \
$\pi_0\in\|F(n_0,1,c)\|$ \ since \ $n_0=n_1$, $\pi_0=\pi_1$.\\
Therefore \ $\kk_{\pi_0}\force\neg F(n_0,0,c)$ \ and \ $\neg F(n_0,1,c)$. Moreover, we have obviously \ $\zeta_0\force0\ne1$, since $\|0\ne1\|=\vide$. Therefore, we have \
$\xi'\star\kk_{\pi_0}\ps\kk_{\pi_0}\ps\zeta_0\ps\pi_0\in\bbot$, which is a
contradiction.
\qed

\smallskip\noindent
Theorems~\ref{pas_surj} and~\ref{B_non_den} show that $\gd$ is infinite and not
equipotent with $\gd\fois\gd$, thus not well orderable. Since $\gd$ is equipotent with a
subset of ${\mathcal P}(\wt{\NN})$ (theorem~\ref{gimel2_R}), we have shown that
${\mathcal P}(\wt{\NN})$ is not well orderable, with the hypothesis of theorem~\ref{pas_surj}.\\
More precisely, by corollary~\ref{gimel_n_R}, we know that
$\gimel\mb{n}$ is equipotent with a subset of ${\mathcal P}(\wt{\NN})$ for each integer $n$.
Therefore, we have~:

\begin{theorem}\label{X_m_fois_X_n}
With the hypothesis of theorem~\ref{pas_surj}, the following formula is realized~:\\
``~There exists a sequence ${\mathcal X}_n$ of infinite subsets of \ ${\mathcal P}(\wt{\NN})$ \ such that, for every integers $m,n\ge2$~:\\
$\bullet$~~there is an injection from ${\mathcal X}_n$ into ${\mathcal X}_{n+1}$~;\\
$\bullet$~~there is no surjection from ${\mathcal X}_n$ onto ${\mathcal X}_{n+1}$~;\\
$\bullet$~~${\mathcal X}_m\fois{\mathcal X}_n$ and ${\mathcal X}_{mn}$ are equipotent~''.\qed
\end{theorem}\noindent
For each integer $n\ge2$, the set $\mb{n}=\{0,1,\ldots,n-1\}$ is a ring~: the ring of integers modulo $n$~; the Boolean algebra $\{0,1\}$ is a set of idempotents in this ring. These ring operations extend to the realizability model,
giving a ring structure on $\gimel\mb{n}$, and $\gd$ is a set of idempotents in~$\gimel\mb{n}$.\\
For each $a\eps\gd$, the equation \ $ax=x$ defines an ideal in $\gimel\mb{n}$, which we denote as $a\gimel\mb{n}$.\\
The application $x\longmapsto ax$ is a retraction from $\gimel\mb{n}$ onto $a\gimel\mb{n}$.

\begin{proposition}\label{aB_mn}
The following formulas are realized in ${\mathcal N}$~:\\
\emph{\phantom ii)}~~$\pt n^{\gn}\pt a^{\gd}($the application \ $x\longmapsto(ax,(1-a)x)$ \ is a bijection\\
\hspace*{\fill}from $\gimel\mb{n}$ onto $a\gimel\mb{n}\,\fois(1-a)\gimel\mb{n})$.\\
\emph{ii)}~~$\pt m^{\gn}\pt n^{\gn}\pt a^{\gd}($the application \ $(x,y)\longmapsto my+x$ is a bijection\\
\hspace*{\fill}from $a\gimel\mb{m}\fois a\gimel\mb{n}$ onto $a\gimel(\mb{mn}))$.
\end{proposition}

\proof\hfill\\
\phantom ii)~~Trivial~: the inverse is \ $(y,y')\longmapsto y+y'$.\\
ii)~~By theorem~\ref{B_mn}, this application is injective~; clearly, it sends \
$a\gimel\mb{m}\fois a\gimel\mb{n}$ \ into \ $a\gimel(\mb{mn})$. Conversely, if \ $z\eps a\gimel(\mb{mn})$,
then there exists $x\eps\gimel\mb{m}$ and $y\eps\gimel\mb{n}$ such that $z=my+x$~;\\
thus, we have $z=az=may+ax$ with \ $ax\eps a\gimel\mb{m}$ \ and \ $ay\eps a\,\gimel\mb{n}$.
\qed

\begin{theorem}\label{non_surj_aB}
We suppose that, for each \ $\alpha\in\Lbd,\,\pi\in\Pi$, and every distinct \
$\zeta_0,\zeta_1,\zeta_2\in\Lbd$, \ we have \ $\kk_\pi\alpha\zeta_0\force\bot$ \ or \ $\kk_\pi\alpha\zeta_1\force\bot$ \ or \ $\kk_\pi\alpha\zeta_2\force\bot$.\\
Then, for each formula $F(x,y,z)$ with three free variables, we have~:

\smallskip\noindent
$\theta\force\pt z\pt m^{\gn}\pt n^{\gn}\pt a^{\gd}[(2m\ppt n)=1\mapto\\
\hspace*{\fill}(a\ne0,\pt x\pt y\pt y'(F(x,y,z),F(x,y',z),y\ne y'\to\bot),
\pt y^{\gimel\mb{n}}\ex x^{\gimel\mb{m}}F(x,ay,z)\to\bot)]$

\smallskip\noindent
with \ $\theta=\lbd a\lbd x\lbd y(\Ccc)\lbd k(y)\lbd z(xzz)(k)az$.
\end{theorem}\noindent
{\small{\bfseries Remark.} This formula means that, if $n>2m,\,a\eps\gd,a\ne0$, then
there is no surjection from $\gimel\mb{m}$ onto $a\gimel\mb{n}$~: it suffices to take \
$F(x,y,z)\equiv{<}x,y{>}\eps z$.}

\smallskip\proof
Reasoning by contradiction, let us consider $m,n\in\NN$ with $n>2m$, $a\in\{0,1\}$, an individual $c$, three terms $\alpha,\xi,\eta\in\Lbd$ and $\pi\in\Pi$ such that~:\\
$\theta\star\alpha\ps\xi\ps\eta\ps\pi\notin\bbot$~, \ $\alpha\force a\ne0$, \
$\xi\force\pt x\pt y\pt y'(F(x,y,c),F(x,y',c),y\ne y'\to\bot)$,\\
$\eta\force\pt y^{\gimel\mb{n}}\neg\pt x^{\gimel\mb{m}}\neg F(x,ay,c)$.\\
We have \ $\theta\star\alpha\ps\xi\ps\eta\ps\pi\succ\eta\star\theta'\ps\pi$ \ and
therefore \ $\eta\star\theta'\ps\pi\notin\bbot$ \ with $\theta'=\lbd z(\xi zz)(\kk_\pi)\alpha z$.\\
It follows that, for every \ $y\in\{0,\ldots,n-1\}$, \ we have \
$\theta'\nforce\pt x^{\gimel\mb{m}}\neg F(x,ay,c)$.\\
Thus, there exist two functions \  $y\longmapsto x_y$ (resp. \ $y\longmapsto\zeta_y$) \
from \ $\{0,\ldots,n-1\}$ \ into \ $\{0,\ldots,m-1\}$ \
(resp. into $\Lbd$), such that \ $\zeta_y\force F(x_y,ay,c)$ \ and \ $\theta'\star\zeta_y\ps\varpi_y\notin\bbot$ (for some suitable stacks $\varpi_y)$.\\
Now, we have \ $\theta'\star\zeta_y\ps\varpi_y\succ\xi\star\zeta_y\ps\zeta_y\ps\kappa_y\ps\varpi_y$ \ with \ $\kappa_y=\kk_\pi\alpha\zeta_y$~; \ therefore, we have~:\\  $\xi\star\zeta_y\ps\zeta_y\ps\kappa_y\ps\varpi_y\notin\bbot$ \ for each $y\in\{0,\ldots,n-1\}$.\\
By hypothesis on $\xi$ (with \ $y=y')$, it follows that \ $\kappa_y\nforce\bot$ \ for every $y<n$.\\
It follows first that \ $\alpha\nforce\bot$ \ and therefore, we have \ $a=1$~; thus \ $\zeta_y\force F(x_y,y,c)$.\\
Moreover, since $n>2m$, there exist \ $y_0,y_1,y_2<n$ distinct, such that \ $x_{y_0}=x_{y_1}=x_{y_2}$.\\
But, following the hypothesis of the theorem, the terms
$\zeta_{y_0},\zeta_{y_1},\zeta_{y_2}$ cannot be distinct, because \ $\kappa_{y_0},\kappa_{y_1},\kappa_{y_2}\nforce\bot$. Therefore we have, for instance, \
$\zeta_{y_0}=\zeta_{y_1}$~; then, we apply the hypothesis on $\xi$ with $y=y_0,y'=y_1$, which gives
$\xi\star\zeta_{y_0}\ps\zeta_{y_1}\ps\kappa\ps\varpi\in\bbot$ \ for every $\kappa\in\Lbd$ and $\varpi\in\Pi$. But it follows that \ $\xi\star\zeta_{y_0}\ps\zeta_{y_0}\ps\kappa_{y_0}\ps\varpi_{y_0}\in\bbot$ \
which is a contradiction.
\qed

\begin{corollary}\label{B_n_aB_n+1}
With the hypothesis of theorem~\ref{non_surj_aB}, the following formulas are realized~:\\
\emph{\phantom{ii}i)}~~$\pt n^{\wt{\NN}}\,\pt a^{\gd}(a\ne0\;\to$ \ there is no surjection from \ $\gimel\mb{n}$ \
onto \ $a\gimel(\mb{n+1}))$.\\
\emph{\phantom iii)}~~$\pt n^{\wt{\NN}}\,\pt a^{\gd}\pt b^{\gd}(a\et b=0,b\ne0\;\to$ \ there is no surjection from \ $a\gimel\mb{n}$ \
onto \ $b\gd)$.\\
\emph{iii)}~~$\pt n^{\wt{\NN}}\,\pt a^{\gd}\pt b^{\gd}(a\et b=a,a\ne b\;\to$ \ there is no surjection from \ $a\gimel\mb{n}$
onto \ $b\gd)$.
\end{corollary}

\proof\hfill\\
\phantom{ii}i)~Suppose that there is a surjection from $\gimel\mb{n}$ onto $a\gimel(\mb{n+1})$. Then,
by the recurrence scheme (theorem~\ref{rec_eps}(ii)), we see that, for each integer $k\in\wt{\NN}$, there exists a surjection from \ $(\gimel\mb{n})^k$ \ onto \ $(a\gimel(\mb{n+1}))^k$~; and, by proposition~\ref{aB_mn}(ii) and the recurrence scheme,
it follows that there is a surjection from \ $\gimel(\mb{n^k})$ \ onto \ $a\gimel(\mb{(n+1)^k})$.\\
But, for $k>n$, we have $(n+1)^k>2n^k$ and this contradicts theorem~\ref{non_surj_aB}.\\
\phantom iii)~~Since $a\et b=0$, the rings $(a+b)\gimel\mb{n}$ and $a\gimel\mb{n}\times b\gimel\mb{n}$
are isomorphic. Reasoning by contradiction, there would exist a surjection from \ $(a+b)\gimel\mb{n}$ \
onto \ $b\gd\fois b\gimel\mb{n}$, thus also onto \ $b\,\gimel(\mb{2n})$ (proposition~\ref{aB_mn}(ii)),
thus a surjection from \ $\gimel\mb{n}$ onto \ $b\gimel(\mb{2n})$, which contradicts~(i).\\
iii)~~Otherwise, there would exist a surjection from $a\gimel\mb{n}$ onto $(b-a)\gd$, which contradicts~(ii).
\qed

\smallskip\noindent
{\bfseries Applications.}\ \\
\phantom ii)~~By DC, since $\gd$ is atomless, there exists in $\gd$ a strictly decreasing sequence. Hence, by  
corollary~\ref{B_n_aB_n+1}(iii) and theorem~\ref{gimel2_R}, there exists a sequence of infinite subsets of
${\mathcal P}(\wt{\NN})$, the ``cardinals'' of which are strictly decreasing.\\
ii)~~Applying corollary~\ref{B_n_aB_n+1}(ii) with $n=2$, we see that there exist two subsets of
${\mathcal P}(\wt{\NN})$ the ``cardinals'' of which are incomparable~; which means that there is no surjection of
one of them onto the other.

\smallskip\noindent
More precisely, let ${\mathcal B}$ be the image of $\gd$ by the injection in ${\mathcal P}(\wt{\NN})$ given by
theorem~\ref{gimel2_R}~; then we have~:

\begin{theorem}\label{subsets_P(N)}
With the hypothesis of theorem~\ref{non_surj_aB}, the following formula is realized
in ${\mathcal N}$~:\\
``There exists a subset ${\mathcal B}$ of ${\mathcal P}(\wt{\NN})$ (the real line of the model ${\mathcal N}$), such that\\
${\mathcal B}$ is an atomless Boolean algebra for the usual order \ $\subseteq$ on
${\mathcal P}(\wt{\NN})$, \ with \ $\vide,\wt{\NN}\in{\mathcal B}$~; \
$a,b\in{\mathcal B}\Fl a\cap b\in{\mathcal B}$.\\
If $a\in{\mathcal B},a\ne\vide$ \ then $a{\mathcal B}$ is infinite and there is no surjection
from \ ${\mathcal B}$ \ onto \ $a{\mathcal B}\fois a{\mathcal B}$\\
(where \ $a{\mathcal B}$ means $\{x\in{\mathcal B};\;x\subseteq a\}$).\\
If $a,b\in{\mathcal B},a,b\ne\vide$ and $a\cap b=\vide$, then there is no surjection
from \ $a{\mathcal B}$ \ onto \ $b{\mathcal B}$ (the ``cardinals'' of \ $a{\mathcal B},b{\mathcal B}$ \ are incomparable).\\
If $a,b\in{\mathcal B},a\subseteq b$ and $a\ne b$, then there is no surjection from
$a{\mathcal B}$ onto $b{\mathcal B}$ (the ``cardinal'' of \ $a{\mathcal B}$ \ is strictly
less than the ``cardinal'' of \ $b{\mathcal B}$)''.\qed
\end{theorem}\noindent
In other words, for $a,b\in{\mathcal B}$, we have~: $a\subseteq b$ $\Dbfl$ there exists a surjection from $b{\mathcal B}$ onto~$a{\mathcal B}$. The order, in the atomless Boolean algebra ${\mathcal B}$, is the order on the
``cardinals'' of its initial segments.

\subsection*{The model of threads}
This model is the canonical instance of a non trivial coherent realizability model. It is defined as follows~:

\smallskip\noindent
Let \ $n\longmapsto\pi_n$ \ be an enumeration of the \emph{stack constants} and let \
$n\longmapsto\theta_n$ \ be a recursive enumeration of the \emph{proof-like terms}.
For each $n\in\NN$, the \emph{thread with number $n$} is the set of processes which
appear during the execution of the process \ $\theta_n\star\pi_n$. In other words,
it is the~set of all processes \ $\xi\star\pi$ \ such that \ $\theta_n\star\pi_n\succ\xi\star\pi$.\\
Note that every term which appears in the $n$-th thread contains the only stack constant~$\pi_n$.

\smallskip\noindent
We define $\bbot^c$ (the complement of $\bbot$) as the union of all threads. Thus, a process
$\xi\star\pi$ is in~$\bbot^c$ \ iff \ $(\ex n\in\NN)\;\theta_n\star\pi_n\succ\xi\star\pi$.\\
Therefore, we have \ $\xi\star\pi\in\bbot$ iff the process $\xi\star\pi$ never appears
in any thread.\\
For every term $\xi$, we have \ $\xi\force\bot$ \ iff \ $\xi$ never appears in head position in any thread.\\
If \ $\xi$ is a proof-like term, we have $\xi=\theta_n$ for some integer $n$, and
therefore $\xi\star\pi_n\notin\bbot$, by definition of $\bbot$. It follows that
\emph{the model of threads is coherent}.

\smallskip\noindent
If $\xi\in\Lbd$, $\xi\nforce\bot$ then $\xi$ appears in head position in at least
one thread. This thread is unique, unless $\xi$ is a proof-like term, because it is determined by the number of any stack constant which appears in $\xi$.

\begin{theorem}
The hypothesis of theorems~\ref{B_non_triviale},~\ref{sans_atome},~\ref{pas_surj} and~\ref{non_surj_aB} are satisfied in the model of threads.
\end{theorem}

\proof
The hypothesis of theorems~\ref{pas_surj} and~\ref{B_non_triviale} are trivially satisfied
if we take~:\\
$\omega=(\lbd x\,xx)\lbd x\,xx$, \ $\omega_0=(\omega)\ul{0}$, \ and \ ${\omega_1=(\omega)\ul{1}}$.\\
Moreover, the hypothesis of theorem~\ref{non_surj_aB} is obviously stronger than the
hypothesis of theorem~\ref{sans_atome}.

\smallskip\noindent
We check the hypothesis of theorem~\ref{non_surj_aB} by contradiction~:\\
Suppose that $\kk_\pi\alpha\zeta_0\nforce\bot$, $\kk_\pi\alpha\zeta_1\nforce\bot$ and $\kk_\pi\alpha\zeta_2\nforce\bot$.
Therefore, these three terms appear in head position, and moreover in the same thread~: indeed, since they contain the stack $\pi$, this thread has the same number as the stack 
constant of~$\pi$.\\
Let us consider their first appearance in head position, for instance with the order $0,1,2$.\\
Therefore we have, in this thread~: \
$\kk_\pi\alpha\zeta_0\star\rho_0\succ\alpha\star\pi\succ\cdots\succ\kk_\pi\alpha\zeta_1\star\rho_1
\succ\alpha\star\pi\succ\cdots$\\
But, at the second appearance of $\alpha\star\pi$, the thread enters into a loop, and the
term \ $\kk_\pi\alpha\zeta_2$ can never arrive in head position, since \ $\zeta_1\ne\zeta_2$.
\qed


\begin{thebibliography}{99}
\bibitem{berardi}S. Berardi, M. Bezem, T. Coquand. {\em On the computational content of the axiom of choice.}
J. Symb. Log. 63 (1998), p. 600-622.

\bibitem{curry}H.B. Curry, R. Feys. {\em Combinatory Logic.} North-Holland (1958).

\bibitem{easton}W. Easton. {\em Powers of regular cardinals.} Ann. Math. Logic 1 (1970), p. 139-178.

\bibitem{fri1} H. Friedman. {\em The consistency of classical set theory relative
to a set theory with intuitionistic logic.} Journal of Symb. Logic, 38 (2)  (1973) p. 315-319.

\bibitem{fri2}  H. Friedman. {\em Classically and intuitionistically provably
recursive functions.}\\
In: Higher set theory. Springer Lect. Notes in Math. 669 (1977) p. 21-27.

\bibitem{girard}J.Y. Girard. {\em Une extension de l'interprétation fonctionnelle de Gödel à l'analyse.}\\
Proc. 2nd Scand. Log. Symp. (North-Holland) (1971) p. 63-92.

\bibitem{griffin}T. Griffin. {\em A formul\ae -as-type notion of control.}\\
Conf. record 17th A.C.M. Symp. on Principles of Progr. Languages (1990).

\bibitem{grigorieff}S. Grigorieff. {\em Combinatorics on ideals and forcing.}\\
Ann. Math. Logic 3(4) (1971), p. 363-394.

\bibitem{howard}W. Howard. {\em The formulas--as--types notion of construction.}\\
Essays on combinatory logic, $\lbd$-calculus, and formalism, J.P. Seldin and J.R. Hindley ed.,
Acad. Press (1980) p. 479--490.

\bibitem{hyland}J. M. E. Hyland. {\em The effective topos.}\\
The L.E.J. Brouwer Centenary Symposium (Noordwijkerhout, 1981), 165--216,\\
Stud. Logic Foundations Math., 110, North-Holland, Amsterdam-New York, 1982.

\bibitem{kreisel1}G. Kreisel. {\em On the interpretation of non-finitist proofs I.}\\
J. Symb. Log. 16 (1951) p. 248-26.

\bibitem{kreisel2}G. Kreisel. {\em On the interpretation of non-finitist proofs II.}\\
J. Symb. Log. 17 (1952), p. 43-58.

\bibitem{krivine1}J.-L. Krivine. {\em Typed lambda-calculus in classical Zermelo-Fraenkel set theory.}\\
Arch. Math. Log., 40, 3, p. 189-205 (2001).\\
http://www.pps.jussieu.fr/\verb?~?krivine/articles/zf\_epsi.pdf

\bibitem{krivine2}J.-L. Krivine. {\em Dependent choice, `quote' and the clock.}\\
Th. Comp. Sc., 308, p. 259-276 (2003).\\
http://hal.archives-ouvertes.fr/hal-00154478\\
Updated version at~:\\
http://www.pps.jussieu.fr/\verb?~?krivine/articles/quote.pdf

\bibitem{krivine3}J.-L. Krivine. {\em Realizability in classical logic.}\\
In {\em Interactive models of computation and program behaviour.}\\
Panoramas et synthèses, Société Mathématique de France, 27, p. 197-229 (2009).\\
http://hal.archives-ouvertes.fr/hal-00154500\\
Updated version at~:\\
http://www.pps.jussieu.fr/\verb?~?krivine/articles/Luminy04.pdf

\bibitem{krivine4}J.-L. Krivine. {\em Realizability : a machine for Analysis and set theory.}\\
Geocal'06 (febr. 2006 - Marseille); Mathlogaps'07 (june 2007 - Aussois).\\
http://cel.archives-ouvertes.fr/cel-00154509\\
Updated version at~:\\
http://www.pps.jussieu.fr/\verb?~?krivine/articles/Mathlog07.pdf

\bibitem{krivine5}J.-L. Krivine. {\em Structures de réalisabilité, RAM et ultrafiltre sur $\NN$.} (2008)\\
http://hal.archives-ouvertes.fr/hal-00321410\\
Updated version at~:\\
http://www.pps.jussieu.fr/\verb?~?krivine/articles/Ultrafiltre.pdf

\bibitem{krivine6}J.-L. Krivine. {\em Realizability algebras~: a program to well order $\mathbb{R}$.}\\
Logical Methods in Computer Science, vol. 7, 3:02, p. 1-47 (2011)\\
http://hal.archives-ouvertes.fr/hal-00483232\\
Updated version at~:\\
http://www.pps.jussieu.fr/\verb?~?krivine/articles/Well\_order.pdf

\end{thebibliography}
\end{document}